\documentclass[11pt]{article}
\usepackage{geometry}                
\usepackage{amsmath,amssymb}
\usepackage{amsthm}
\usepackage{hyperref}
\usepackage{graphicx,epsfig,color}
\usepackage{booktabs,bm,multirow}
\usepackage{cases}
\usepackage[caption=false]{subfig}

\newtheorem{theorem}{Theorem}[section]

\newtheorem{remark}[theorem]{Remark}
\newtheorem{definition}[theorem]{Definition}
\newtheorem{lemma}[theorem]{Lemma}

\def\mcalp{\mathcal P}
\def\mcalx{\mathcal X}
\DeclareMathOperator*{\argmin}{argmin}
\def\mbbe{\mathbb E}
\def\mbbr{\mathbb R}
\def\p{\partial}
\def\la{\langle} \def\ra{\rangle}
\def\l{\left} 
\def\r{\right} 
\def\mbf{\mathbf}
\def\rmf{{\rm F}}
\def\wh{\widehat}
\def\rank{{\rm rank}} 
\def\ran{{\rm range}} 
\def\nul{{\rm null}}
\def\sgn{{\rm sgn}} 
\def\dsp{\displaystyle} 
\def\nn{\nonumber}
\newcommand{\beq}{\begin{equation}}\newcommand{\eeq}{\end{equation}}
\newcommand{\bem}{\begin{bmatrix}}\newcommand{\eem}{\end{bmatrix}}

\title{Regularized randomized iterative algorithms for factorized linear systems}
\author{Kui Du\thanks{School of Mathematical Sciences and Fujian Provincial Key Laboratory of Mathematical Modeling and High Performance Scientific Computing, Xiamen University, Xiamen 361005, China (kuidu@xmu.edu.cn).}} 
\date{}

\begin{document}
\maketitle

\begin{abstract} Randomized iterative algorithms for solving the factorized linear system, $\mbf A\mbf B\mbf x=\mbf b$ with $\mbf A\in\mbbr^{m\times \ell}$, $\mbf B\in\mbbr^{\ell\times n}$, and $\mbf b\in\mbbr^m$, have recently been proposed. They take advantage of the factorized form and avoid forming the matrix $\mbf C=\mbf A\mbf B$ explicitly. However, they can only find the minimum norm (least squares) solution. In contrast, the regularized randomized Kaczmarz (RRK) algorithm can find solutions with certain structures from consistent linear systems. In this work, by combining the randomized Kaczmarz algorithm or the randomized Gauss--Seidel algorithm with the RRK algorithm, we propose two new regularized randomized iterative algorithms to find (least squares) solutions with certain structures of $\mbf A\mbf B\mbf x=\mbf b$. We prove linear convergence of the new algorithms. Computed examples are given to illustrate that the new algorithms can find sparse (least squares) solutions of $\mbf A\mbf B\mbf x=\mbf b$ and can be better than the existing randomized iterative algorithms for the corresponding full linear system $\mbf C\mbf x=\mbf b$ with $\mbf C=\mbf A\mbf B$.  
\vspace{.5mm} 



\end{abstract}

\flushbottom
\section{Introduction}  

We are interested in solving the following large-scale factorized linear system \beq\label{fls} \bf ABx=b,\eeq where $$\mbf A\in\mbbr^{m\times\ell},\quad \mbf B\in\mbbr^{\ell\times n}, \quad \rank(\mbf A)=\rank(\mbf B)=\ell,\quad \mbf b\in\mbbr^m.$$ This kind of system either arises naturally in many applications or may be imposed to save space required for storage when working with large low-rank datasets; see, for example,  \cite{ma2018itera} and the references therein. For the problems of interest, $m$ and $n$ may be so large that existing methods that require all-at-once access to the matrix $\mbf C=\mbf A\mbf B$ are not feasible. Instead, we consider randomized or sampling methods where only ``blocks'' of the matrices $\mbf A$ and $\mbf B$ are required at a given time. 

By introducing an auxiliary vector $\mbf y$, the system (\ref{fls}) can be written as two individual subsystems $\mbf A\mbf y=\mbf b$ (possibly inconsistent) and $\mbf B\mbf x=\mbf y$ (always consistent). Then one can find a solution of (\ref{fls}) by solving each subsystem separately. However, if iterative methods are used, it is usually unclear when the iterates of the first subsystem should be stopped (if terminated prematurely, the error may propagate through iterates when solving the second subsystem). Randomized iterative algorithms that can address this issue and take advantage of the factorized form have recently been proposed. For example, by intertwining the randomized Kaczmarz (RK) algorithm \cite{strohmer2009rando} or the randomized extended Kacamzrz (REK) algorithm \cite{zouzias2013rando} for the subsystem $\mbf A\mbf y=\mbf b$ with the RK algorithm for the subsystem $\mbf B\mbf x=\mbf y$ in an alternating fashion, Ma et al. \cite{ma2018itera} proposed the RK-RK algorithm and the REK-RK algorithm for the system (\ref{fls}). Similarly, by intertwining the randomized Gauss--Seidel (RGS) algorithm \cite{leventhal2010rando} with the RK algorithm, Zhao et al. \cite{zhao2023rando} proposed the RGS-RK algorithm for the system (\ref{fls}). The RK-RK algorithm converges linearly to the unique minimum norm solution if the system (\ref{fls}) is consistent. The REK-RK and RGS-RK algorithms converge linearly to the unique minimum norm least squares solution if the system (\ref{fls}) is inconsistent. Note that these algorithms avoid forming the matrix $\mbf C=\mbf A\mbf B$ explicitly.

In many applications, the matrix $\mbf B$ in (\ref{fls}) acts as a frame or redundant dictionary, and the system (\ref{fls}) may have a sparse solution (usually not the minimum norm one). If the matrix $\mbf C=\mbf A\mbf B$ is explicitly given, the well known randomized sparse Kaczmarz (RSK) algorithm \cite{lorenz2014spars,petra2015rando,schopfer2019linea,chen2021regul} can be used to find a sparse solution for the consistent case, and the recently proposed generalized extended randomized Kaczmarz (GERK-(a,d)) algorithm \cite{schopfer2022exten} can be used to find a sparse least squares solution for the inconsistent case. It was shown in \cite{schopfer2019linea} that for the consistent case the RSK algorithm converges linearly to the unique solution of the regularized basis pursuit problem $$\mbox{minimize } \frac{1}{2}\|\mbf x\|_2^2+\lambda\|\mbf x\|_1,\quad \mbox{s.t.}\quad \mbf C \mbf x=\mbf b.$$ It was shown in \cite{schopfer2022exten} that the GERK-(a,d) algorithm converges linearly to the unique solution of the combined optimization problem \beq\label{cls} \mbox{minimize } \frac{1}{2}\|\mbf x\|_2^2+\lambda\|\mbf x\|_1,\quad \mbox{s.t.}\quad \mbf x\in\argmin_{\mbf z\in\mbbr^n} \|\mbf b-\mbf C \mbf z\|_2.\eeq  Actually, the GERK-(a,d) algorithm for (\ref{cls}) can be viewed as an RK-RSK approach, which combines the RK algorithm (with initial iterate $\mbf z^{(0)}=\mbf b$) for $\mbf C^\top\mbf z=\mbf 0$ and the RSK algorithm for  $\mbf C\mbf x = \mbf C\mbf C^\dag\mbf b$. Here $(\cdot)^\dag$ denotes the Moore--Penrose pseudoinverse. 

In this paper, we adopt the idea used by the GERK-(a,d) algorithm and propose two regularized randomized iterative algorithms for solving the following combined optimization problem: \beq\label{mp} \mbox{minimize } f(\mbf x),\quad \mbox{s.t.}\quad \mbf x\in\argmin_{\mbf z\in\mbbr^n} \|\mbf b-\mbf A\mbf B \mbf z\|_2,\eeq where the objective function $f$ is strongly convex, and $\mbf A$, $\mbf B$, and $\mbf b$ are as given in (\ref{fls}). The proposed algorithms are ``regularized'' since the objective function usually contains regularization terms  for promoting certain structures of the underlying solution. For example, the strongly convex function $\dsp f(\mbf x)=\frac{1}{2}\|\mbf x\|_2^2+\lambda\|\mbf x\|_1$ with $\lambda>0$ can be used to promote sparsity. Specifically, our proposed algorithms intertwine the RK algorithm or the RGS algorithm for the subsystem $\mbf A\mbf y= \mbf b$ with the regularized randomized Kaczmarz (RRK) algorithm \cite{lorenz2014spars,petra2015rando,schopfer2019linea,chen2021regul} for the linear equality constrained minimization problem \beq\label{mp2}\mbox{minimize } f(\mbf x),\quad \mbox{s.t.}\quad \mbf B\mbf x=\mbf y.\eeq They avoid forming the matrix $\mbf C=\mbf A\mbf B$ explicitly,  only require a row or column of $\mbf  A$ and a row of $\mbf B$ at each step, and become the RK-RK algorithm and the RGS-RK algorithm if we set $f(\mbf x)=\dsp \frac{1}{2}\|\mbf x\|_2^2$ in (\ref{mp2}).

The rest of this paper is organized as follows. In section 2, we provide clarification of notation, a brief review of fundamental concepts and results in convex optimization, and the basic results for the RK, RGS, and RRK algorithms. In section 3, we  introduce the proposed algorithms and prove their linear convergence in expectation to the unique solution of $(\ref{mp})$. We also construct an extension to the RGS algorithm that parallels the GERK-(a,d) algorithm, which converges linearly to the unique solution of (\ref{cls}). In section 4, some numerical examples are performed to demonstrate the theoretical results and the effectiveness of our algorithms. Finally, we give some concluding remarks and possible future work in section 5.

\section{Preliminaries} 

\subsection{Notation} For any random variable $\bm\xi$, we use $\mbbe\bem\bm\xi\eem$ to denote the expectation of $\bm\xi$. For an integer $n\geq 1$, let $[n]$ denote the set $\{1,2,\cdots,n\}$. For any column vector $\mbf b\in\mbbr^m$, we use $b_i$, $\mbf b^\top$, $\|\mbf b\|_1$, and $\|\mbf b\|_2$ to denote the $i$th component, the transpose, the $\ell_1$ norm, and the Euclidean norm of $\mbf b$, respectively. We use $\mbf I$ to denote the identity matrix whose order is clear from the context. For any matrix $\mbf A\in\mbbr^{m\times\ell}$, we use $\mbf A_{i,:}$, $\mbf A_{:,j}$, $\mbf A^\top$, $\mbf A^\dag$, $\|\mbf A\|_2$, $\|\mbf A\|_\rmf$, $\ran(\mbf A)$, $\nul(\mbf A)$, $\rank(\mbf A)$, and $\sigma_{\rm min}(\mbf A)$ to denote the $i$th row, the $j$th column, the transpose, the Moore--Penrose pseudoinverse, the 2-norm, the Frobenius norm, the column space, the null space, the rank, and the minimum nonzero singular value of $\mbf A$, respectively. For any $\mbf x,\mbf y\in\mbbr^m$, we use $\la\mbf x,\mbf y\ra$ to denote the standard inner product, i.e., $\la\mbf x,\mbf y\ra=\mbf x^\top\mbf y$. Define the soft shrinkage function $S_\lambda(\mbf x)$ component-wise as $$(S_\lambda(\mbf x))_i= \max\{|x_i|-\lambda, 0\}\sgn(x_i),$$ where $\sgn(\cdot)$ is the sign function.

\subsection{Convex optimization basics} To make the paper self-contained, we review basic definitions and properties about convex functions defined on $\mbbr^n$ in this subsection. We refer the reader to \cite{rockafellar1970conve,beck2017first} for more definitions and properties. 

\begin{definition}[subdifferential]
For a function $f:\mbbr^n\rightarrow\mbbr$, its subdifferential at $\mbf x\in\mbbr^n$ is defined as $$\p f(\mbf x):=\{\mbf z\in\mbbr^n\mid f(\mbf y)\geq f(\mbf x)+\la\mbf z,\mbf y-\mbf x\ra, \quad  \forall\ \mbf y\in\mbbr^n\}.$$	
\end{definition}

\begin{definition}[$\gamma$-strong convexity]
A function $f:\mbbr^n\rightarrow\mbbr$ is called $\gamma$-strongly convex for a given $\gamma>0$ if the following inequality holds for all $\mbf x, \mbf y\in\mbbr^n$ and $\mbf z\in\p f(\mbf x)$: $$f(\mbf y)\geq f(\mbf x)+\la\mbf z,\mbf y-\mbf x\ra+\frac{\gamma}{2}\|\mbf y-\mbf x\|_2^2.$$	
\end{definition}
As an example, the function $f(\mbf x)=\dsp\frac{1}{2}\|\mbf x\|_2^2+\lambda\|\mbf x\|_1$ with $\lambda\geq 0$ is $1$-strongly convex.

\begin{definition}[conjugate function]
The conjugate function of $f:\mbbr^n\rightarrow\mbbr$ at $\mbf x\in\mbbr^n$ is defined as $$f^*(\mbf x):=\sup_{\mbf y\in\mbbr^n}\{\la\mbf x,\mbf y\ra-f(\mbf y)\}.$$	
\end{definition} 
If $f(\mbf x)$ is $\gamma$-strongly convex, then the conjugate function $f^*(\mbf x)$ is differentiable and for all $\mbf x,\mbf y\in\mbbr^n$, the following inequality holds: 
 \beq\label{sc2}f^*(\mbf y)\leq f^*(\mbf x)+\la\nabla f^*(\mbf x),\mbf y-\mbf x\ra+\frac{1}{2\gamma}\|\mbf y-\mbf x\|_2^2.\eeq For a strongly convex function $f(\mbf x)$, it can be shown that \cite{rockafellar1970conve,beck2017first} \beq\label{xgz}\mbf z\in\p f(\mbf x) \ \Leftrightarrow\ \mbf x=\nabla f^*(\mbf z).\eeq

\begin{definition}[Bregman distance] For a convex function $f:\mbbr^n\rightarrow\mbbr$, the Bregman distance between $\mbf x$ and $\mbf y$ with respect to $f$ and $\mbf z\in \p f(\mbf x)$ is defined as $$D_{f,\mbf z}(\mbf x,\mbf y):=f(\mbf y)-f(\mbf x)-\la\mbf z,\mbf y-\mbf x\ra.$$ 
\end{definition}
It follows from $\mbf z\in \p f(\mbf x)$ that $\la\mbf z,\mbf x\ra=f(\mbf x)+f^*(\mbf z)$ (see \cite[Theorem 4.20]{beck2017first}). Then it holds that \beq\label{dfz} D_{f,\mbf z}(\mbf x,\mbf y)=f(\mbf y)+f^*(\mbf z)-\la\mbf z,\mbf y\ra.\eeq If $f$ is $\gamma$-strongly convex, then it holds that \beq\label{gamma} D_{f,\mbf z}(\mbf x,\mbf y)\geq \frac{\gamma}{2}\|\mbf x-\mbf y\|_2^2.\eeq

\begin{definition}[restricted strong convexity \cite{lai2013augme,schopfer2016linea}]
Let $f:\mbbr^n\rightarrow\mbbr$ be convex differentiable with a nonempty minimizer set $\mcalx_f$. The function $f$ is called restricted strongly convex on $\mbbr^n$ with a constant $\mu>0$ if it satisfies for all $\mbf x\in \mbbr^n$ the inequality $$\la \nabla f(\mcalp_f(\mbf x))-\nabla f(\mbf x), \mcalp_f(\mbf x)-\mbf x\ra\geq \mu\|\mcalp_f(\mbf x)-\mbf x\|_2^2,$$ where $\mcalp_f(\mbf x)$ denotes the orthogonal projection of $\mbf x$ onto $\mcalx_f$. 
\end{definition}

If $\mbf b\in\ran(\mbf A\mbf B)$, then the dual problem of (\ref{mp}) is the unconstrained problem \beq\label{g1}\min_{\mbf y\in\mbbr^m} g(\mbf y):=f^*(\mbf B^\top\mbf A^\top\mbf y)-\la\mbf y,\mbf b\ra.\eeq If $\mbf b\notin\ran(\mbf A\mbf B)$, then the dual problem of (\ref{mp}) is the unconstrained problem \beq\label{g2}\min_{\mbf y\in\mbbr^n} g(\mbf y):=f^*(\mbf B^\top\mbf A^\top\mbf A\mbf B\mbf y)-\la\mbf y,\mbf B^\top\mbf A^\top\mbf b\ra.\eeq 

\begin{definition}[strong admissibility]
Let $\mbf A\in\mbbr^{m\times \ell}$, $\mbf B\in\mbbr^{\ell \times n}$ and $\mbf b\in\mbbr^m$ be given. Let $f:\mbbr^n\rightarrow\mbbr$ be strongly convex. The function $f$ is called strongly admissible for $(\mbf A, \mbf B, \mbf b)$ if the function $g(\mbf y)$ defined in {\rm(\ref{g1})} for $\mbf b\in\ran(\mbf A\mbf B)$ or {\rm(\ref{g2})} for $\mbf b\notin\ran(\mbf A\mbf B)$ is restricted strongly convex on $\mbbr^m$ or $\mbbr^n$.
\end{definition}

As an example, the function $\dsp f(\mbf x)=\frac{1}{2}\|\mbf x\|_2^2+\lambda\|\mbf x\|_1$ is strongly admissible for $(\mbf A, \mbf B, \mbf b)$ (see \cite[Example 3.7]{chen2021regul} and \cite[Lemma 4.6]{lai2013augme}). We refer the reader to \cite{schopfer2016linea} for more examples of strongly admissible functions. The following property of strongly admissible functions  is important for our analysis. 

\begin{lemma}\label{lemmanu} Let $\mbf x_\star$ be the solution of {\rm(\ref{mp})}. If $f$ is strongly admissible for $(\mbf A, \mbf B, \mbf b)$, then there exists a constant $\nu>0$ such that \beq\label{nu} D_{f,\mbf z}(\mbf x,\mbf x_\star)\leq \frac{1}{\nu}\|\mbf B(\mbf x-\mbf x_\star)\|_2^2,\eeq for all $\mbf x\in \mbbr^n$ and $\mbf z\in\p f(\mbf x)\cap\ran(\mbf B^\top)$.
\end{lemma}	
\proof
Case (i): $\mbf b\in\ran(\mbf A\mbf B)$. The solution $\mbf x_\star$ of (\ref{mp}) satisfies  $\mbf A\mbf B\mbf x_\star=\mbf b$. By the strong duality, we have $$f(\mbf x_\star)=-\min_{\mbf y \in\mbbr^m}g(\mbf y).$$ Since $\mbf z\in\ran(\mbf B^\top)=\ran(\mbf B^\top\mbf A^\top)$, we can write $\mbf z=\mbf B^\top\mbf A^\top\mbf u$ for some $\mbf u\in\mbbr^m$. Then \begin{align*} D_{f,\mbf z}(\mbf x,\mbf x_\star) &  \stackrel{(\ref{dfz})}{=} f^*(\mbf z)-\la\mbf z,\mbf x_\star\ra+f(\mbf x_\star)\\ &=f^*(\mbf B^\top\mbf A^\top\mbf u) -\la\mbf B^\top\mbf A^\top\mbf u,\mbf x_\star\ra +f(\mbf x_\star)\\ &= f^*(\mbf B^\top\mbf A^\top\mbf u) -\la\mbf u,\mbf A\mbf B\mbf x_\star\ra +f(\mbf x_\star)\\ &= f^*(\mbf B^\top\mbf A^\top\mbf u)-\la\mbf u,\mbf b\ra +f(\mbf x_\star)\\ &= g(\mbf u)-\min_{\mbf y \in\mbbr^m}g(\mbf y).\end{align*} Since $g(\mbf y)$ is restricted strongly convex on $\mbbr^m$, there exists a constant $\mu>0$ such that $$\la \nabla g(\mcalp_g(\mbf u))-\nabla g(\mbf u), \mcalp_g(\mbf u)-\mbf u\ra\geq \mu\|\mcalp_g(\mbf u)-\mbf u\|_2^2.$$ By the optimality condition $\nabla g(\mcalp_g(\mbf u))=0$ and the Cauchy--Schwarz inequality, we get $$\|\nabla g(\mbf u)\|_2\geq\mu\|\mcalp_g(\mbf u)-\mbf u\|_2.$$  The convexity of $g(\mbf y)$ implies \beq\label{convex} g(\mbf u)-g(\mcalp_g(\mbf u))\leq \la\nabla g(\mbf u),\mbf u-\mcalp_g(\mbf u)\ra\leq\|\nabla g(\mbf u)\|_2\|\mcalp_g(\mbf u)-\mbf u\|_2\leq \frac{1}{\mu}\|\nabla g(\mbf u)\|_2^2.\eeq The gradient of $g(\mbf y)$ at $\mbf y=\mbf u$ is \beq\label{grad} \nabla g(\mbf u)=\mbf A\mbf B\nabla f^*(\mbf B^\top\mbf A^\top\mbf u)-\mbf b= \mbf A\mbf B\nabla f^*(\mbf z)-\mbf b  \stackrel{(\ref{xgz})}{=} \mbf A\mbf B\mbf x-\mbf b= \mbf A\mbf B\mbf x-\mbf A\mbf B\mbf x_\star.\eeq Therefore, \begin{align*} D_{f,\mbf z}(\mbf x,\mbf x_\star) &= g(\mbf u)-\min_{\mbf y \in\mbbr^m}g(\mbf y)=g(\mbf u)-g(\mcalp_g(\mbf u))\\ &\stackrel{(\ref{convex})(\ref{grad})}\leq \frac{1}{\mu}\|\mbf A\mbf B(\mbf x-\mbf x_\star)\|_2^2\leq \frac{\|\mbf A\|_2^2}{\mu}\|\mbf B(\mbf x-\mbf x_\star)\|_2^2.\end{align*} 

Case (ii): $\mbf b\notin\ran(\mbf A\mbf B)$. The solution $\mbf x_\star$ of (\ref{mp}) satisfies $\mbf B^\top\mbf A^\top\mbf A\mbf B\mbf x_\star=\mbf B^\top\mbf A^\top\mbf b.$ Note that $$\mbf z\in\ran(\mbf B^\top)=\ran(\mbf B^\top\mbf A^\top)=\ran(\mbf B^\top\mbf A^\top\mbf A\mbf B).$$ By the similar argument as used in Case (i), we have $$ D_{f,\mbf z}(\mbf x,\mbf x_\star) \leq \frac{1}{\mu}\|\mbf B^\top\mbf A^\top\mbf A\mbf B(\mbf x-\mbf x_\star)\|_2^2\leq \frac{\|\mbf B^\top\mbf A^\top\mbf A\|_2^2}{\mu}\|\mbf B(\mbf x-\mbf x_\star)\|_2^2.$$ This completes the proof. 
\endproof

 The constant $\nu$ depends on the matrices $\mbf A$ and $\mbf B$ and the strongly admissible function $f$. For the choice of $\dsp f(\mbf x)=\frac{1}{2}\|\mbf x\|_2^2+\lambda\|\mbf x\|_1$, explicit expressions for $\nu$ can be obtained (see \cite{lai2013augme,schopfer2022exten}). In general, it is hard to quantify $\nu$. 
 
\subsection{Randomized Kaczmarz (RK)} At each step, the RK algorithm \cite{strohmer2009rando} for solving $\mbf A\mbf y = \mbf b$ orthogonally projects the current estimate vector onto the affine hyperplane defined by a randomly chosen row of $\mbf A\mbf y = \mbf b$. See Algorithm 1 for details.

\begin{center}
\begin{tabular*}{160mm}{l}
\toprule {\bf Algorithm 1:} The RK algorithm for solving $\mbf A\mbf y = \mbf b$\\ 
\hline \noalign{\smallskip}
\quad {\bf Input}: $\mbf A\in\mbbr^{m\times \ell}$,  $\mbf b\in\mbbr^m$, and maximum number of  iterations {\tt maxit}.\\\noalign{\smallskip}
\quad {\bf Output}: an approximation of the solution of $\mbf A\mbf y = \mbf b$.
\\ \noalign{\smallskip}
\quad {\bf Initialize}: $\mbf y^{(0)}\in\mbbr^\ell$.\\ \noalign{\smallskip}
\quad {\bf for} $k=1,2,\ldots,$ {\tt maxit} {\bf do}\\ \noalign{\smallskip}
\quad \qquad  Pick $j_k\in[m]$ with probability ${\|\mbf A_{j_k,:}\|^2_2}/{\|\mbf A\|_\rmf^2}$\\  \noalign{\smallskip}
\quad \qquad  Set $\mbf y^{(k)}=\mbf y^{(k-1)}-\dsp\frac{\mbf A_{j_k,:}\mbf y^{(k-1)}-b_{j_k}}{\|\mbf A_{j_k,:}\|_2^2}(\mbf A_{j_k,:})^\top$\\  \noalign{\smallskip}
\quad {\bf end}\\
\bottomrule
\end{tabular*}
\end{center}

If $\mbf b \in \ran(\mbf A)$ (i.e., $\mbf A\mbf y=\mbf b$ is consistent), the $k$th iterate $\mbf y^{(k)}$  in the RK algorithm with arbitrary $\mbf y^{(0)}\in\mbbr^\ell$ satisfies \beq\label{az}\mbbe\bem\|\mbf y^{(k)}-\mbf y^{(0)}_\star\|_2^2\eem\leq\alpha^k\|\mbf y^{(0)}-\mbf y^{(0)}_\star\|_2^2\eeq with $$\alpha=1-\frac{\sigma_{\rm min}^2(\mbf A)}{\|\mbf A\|_\rmf^2},\qquad \mbf y^{(0)}_\star=(\mbf I-\mbf A^\dag\mbf A)\mbf y^{(0)}+\mbf A^\dag\mbf b.$$ This means that the RK algorithm converges linearly to $\mbf y^{(0)}_\star$, the orthogonal projection of the initial guess $\mbf y^{(0)}$ onto the solution set $\{\mbf y\in\mbbr^\ell \mid \mbf A\mbf y=\mbf b\}$.

\subsection{Randomized Gauss--Seidel (RGS)}

At each step, the RGS algorithm \cite{leventhal2010rando} for solving $\min_{\mbf y}\|\mbf b-\mbf A\mbf y\|_2$ updates one component of the current estimate vector by using a randomly chosen column of $\mbf A$. See Algorithm 2 for details.

\begin{center}
\begin{tabular*}{160mm}{l}
\toprule {\bf Algorithm 2:} The RGS algorithm for solving $\min_{\mbf y}\|\mbf b-\mbf A\mbf y\|_2$\\ 
\hline \noalign{\smallskip}
\quad {\bf Input}: $\mbf A\in\mbbr^{m\times \ell}$, $\mbf b\in\mbbr^m$, and maximum number of  iterations {\tt maxit}.\\ \noalign{\smallskip}
\quad {\bf Output}: an approximation of the solution of $\min_{\mbf y}\|\mbf b-\mbf A\mbf y\|_2$.
\\ \noalign{\smallskip}
\quad {\bf Initialize}: $\mbf y^{(0)}\in\mbbr^\ell$.\\ \noalign{\smallskip}
\quad {\bf for} $k=1,2,\ldots,$ {\tt maxit} {\bf do}\\ \noalign{\smallskip}
\quad \qquad Pick $j_k\in[\ell]$ with probability ${\|\mbf A_{:,j_k}\|^2_2}/{\|\mbf A\|_\rmf^2}$\\  \noalign{\smallskip}
\quad \qquad  Set $\dsp\mbf y^{(k)}=\mbf y^{(k-1)}+\frac{(\mbf A_{:,j_k})^\top(\mbf b-\mbf A\mbf y^{(k-1)})}{\|\mbf A_{:,j_k}\|_2^2}\mbf I_{:,j_k}$\\  \noalign{\smallskip}
\quad {\bf end}\\
\bottomrule
\end{tabular*}
\end{center}

For all $\mbf y^{(0)}\in\mbbr^\ell$, the $k$th iterate $\mbf y^{(k)}$ in the RGS algorithm satisfies $$\mbbe\bem\|\mbf A(\mbf y^{(k)}-\mbf A^\dag\mbf b)\|_2^2\eem\leq\alpha^k\|\mbf A(\mbf y^{(0)}-\mbf A^\dag\mbf b)\|_2^2$$ with $$\alpha=1-\frac{\sigma_{\rm min}^2(\mbf A)}{\|\mbf A\|_\rmf^2}.$$ If $\mbf A$ has full column rank, then it follows that \begin{align}\label{rgsbound} \mbbe\bem\|\mbf y^{(k)}-\mbf A^\dag\mbf b\|_2^2\eem &=\mbbe\bem\|\mbf A^\dag\mbf A(\mbf y^{(k)}-\mbf A^\dag\mbf b)\|_2^2\eem \nn \\ & \leq\|\mbf A^\dag\|_2^2\mbbe\bem\|\mbf A(\mbf y^{(k)}-\mbf A^\dag\mbf b)\|_2^2\eem  \leq\alpha^k\|\mbf A^\dag\|_2^2\|\mbf A(\mbf y^{(0)}-\mbf A^\dag\mbf b)\|_2^2.\end{align} This means that the RGS algorithm converges linearly to the unique least squares solution $\mbf A^\dag\mbf b$.

\subsection{Regularized randomized Kaczmarz (RRK)}
Let $f$ be a given $\gamma$-strongly convex function and $\mbf B$ be a given matrix. The following RRK algorithm \cite{lorenz2014spars,petra2015rando,schopfer2019linea,chen2021regul} has been proposed for solving the minimization problem (\ref{mp2}). 

\begin{center}
\begin{tabular*}{160mm}{l}
\toprule {\bf Algorithm 3:} The RRK algorithm for solving $\min_{\mbf x} f(\mbf x)\ \mbox{s.t.}\ \mbf B\mbf x=\mbf y$\\ 
\hline \noalign{\smallskip}
\quad {\bf Input}: $\mbf B\in\mbbr^{\ell\times n}$,  $\mbf y\in\mbbr^\ell$, $\gamma>0$, and maximum number of  iterations {\tt maxit}.\\\noalign{\smallskip}
\quad {\bf Output}: an approximation of the solution of $\min_{\mbf x} f(\mbf x)\ \mbox{s.t.}\ \mbf B\mbf x=\mbf y$.
\\ \noalign{\smallskip}
\quad {\bf Initialize}: $\mbf z^{(0)}\in\ran(\mbf B^\top)$ and $\mbf x^{(0)}=\nabla f^*(\mbf z^{(0)})$.\\ \noalign{\smallskip}
\quad {\bf for} $k=1,2,\ldots,$ {\tt maxit} {\bf do}\\ \noalign{\smallskip}
\quad \qquad  Pick $i_k\in[\ell]$ with probability ${\|\mbf B_{i_k,:}\|^2_2}/{\|\mbf B\|_\rmf^2}$\\  \noalign{\smallskip}
\quad \qquad Set $\mbf z^{(k)} = \mbf z^{(k-1)}-\gamma\dsp \frac{\mbf B_{i_k,:}\mbf x^{(k-1)}-y_{i_k}}{\|\mbf B_{i_k,:}\|^2_2}(\mbf B_{i_k,:})^\top$ \\  \noalign{\smallskip}
\quad \qquad Set $\mbf x^{(k)}=\nabla f^*(\mbf z^{(k)})$ \\  \noalign{\smallskip}
\quad {\bf end}\\
\bottomrule
\end{tabular*}
\end{center}

It was observed (see \cite{petra2015rando,schopfer2019linea}) that the RRK algorithm could be viewed as a random coordinate descent algorithm \cite{nesterov2012effic} applied to the dual objective function $$f^*(\mbf B^\top \mbf u)-\la\mbf u,\mbf y\ra.$$   Specifically,  a negative stochastic gradient step in the random $i_k$th component of $\mbf u$ is given as $$\mbf u^{(k)}=\mbf u^{(k-1)}-\frac{\mbf B_{i_k,:}\nabla f^*(\mbf B^\top\mbf u^{(k-1)})-y_{i_k}}{\|\mbf B_{i_k,:}\|_2^2}\mbf I_{:,i_k}.$$ We can easily recover the RRK algorithm by simply introducing $\mbf z^{(k)} = \mbf B^\top\mbf u^{(k)}$ and using the relation between the primal and dual variables given by $\mbf x= \nabla f^*(\mbf B^\top\mbf u)$. 

Let $\mbf x_\star$ be the unique solution of the minimization problem (\ref{mp2}). If the objective function $f$ is $\gamma$-strongly convex and satisfies, for all $\mbf x\in\mbbr^n$ and $\mbf z\in\p f(\mbf x)\cap\ran(\mbf B^\top)$, $$D_{f,\mbf z}(\mbf x,\mbf x_\star)\leq\frac{1}{\nu_0}\|\mbf B(\mbf x-\mbf x_\star)\|_2^2,$$ then for all $\mbf z^{(0)}\in\ran(\mbf B^\top)$, the sequences $\{\mbf x^{(k)}\}$ and $\{\mbf z^{(k)}\}$ in the RRK algorithm satisfy $$\mbbe\bem D_{f,\mbf z^{(k)}}(\mbf x^{(k)},\mbf x_\star)\eem\leq\beta_0^k D_{f,\mbf z^{(0)}}(\mbf x^{(0)},\mbf x_\star)$$  with $$\beta_0=1-\frac{\gamma\nu_0}{2\|\mbf B\|_\rmf^2}.$$ It follows from (\ref{gamma}) that $$\mbbe\bem\|\mbf x^{(k)}-\mbf x_\star\|_2^2\eem\leq\beta_0^k \frac{2}{\gamma}D_{f,\mbf z^{(0)}}(\mbf x^{(0)},\mbf x_\star).$$ This means that the RRK algorithm converges linearly to the unique solution of the minimization problem (\ref{mp2}). 

For the concrete choices of $\dsp f(\mbf x)=\frac{1}{2}\|\mbf x\|_2^2$ and $\dsp f(\mbf x)=\frac{1}{2}\|\mbf x\|_2^2+\lambda\|\mbf x\|_1$, we have $\nabla f^*(\mbf x)=\mbf x$ and $\nabla f^*(\mbf x)=S_\lambda(\mbf x)$ (see, e.g.,  \cite{beck2017first}), respectively. As a direct consequence, the RRK algorithm becomes the RK algorithm and the RSK algorithm, respectively. 

\section{The proposed algorithms}
Like the algorithms in the works \cite{zouzias2013rando,ma2018itera,du2019tight,zhao2023rando,schopfer2022exten}, our approach for solving the optimization problem (\ref{mp}) combines two randomized iterative algorithms. Specifically, for the consistent case ($\mbf b\in\ran(\mbf A\mbf B)$), we propose using the RK algorithm to solve the subsystem $\mbf A\mbf y=\mbf b$ followed by the RRK algorithm to solve the minimization problem (\ref{mp2}) as shown in Algorithm 4, and call it the RK-RRK algorithm. For the inconsistent case ($\mbf b\notin\ran(\mbf A\mbf B)$), we propose using the RGS algorithm to solve the problem $\min_{\mbf y}\|\mbf b-\mbf A\mbf y\|_2$ followed by the RRK algorithm to solve the minimization problem (\ref{mp2}) as shown in Algorithm 5, and call it the RGS-RRK algorithm.

\subsection{The RK-RRK algorithm for the case $\mbf b\in\ran(\mbf A\mbf B)$}
In this subsection, we analyze the RK-RRK algorithm (Algorithm 4) and prove its linear convergence property. We note that $\mbf y^{(k)}$ in the RK-RRK algorithm is actually the $k$th iterate of the RK algorithm for $\mbf A\mbf y=\mbf b$, whose convergence estimate is given in (\ref{az}). The vectors $\mbf z^{(k)}$ and $\mbf x^{(k)}$ are one-step updates of the RRK algorithm for  $\mbf B\mbf x=\mbf y^{(k)}$ from $\mbf z^{(k-1)}$ and $\mbf x^{(k-1)}$.  The convergence result of the RK-RRK algorithm is given in Theorem \ref{main}. The proof uses the same idea as that used in \cite{du2019tight,du2020rando}, but is slightly more complicated.

\begin{center}
\begin{tabular*}{160mm}{l}
\toprule {\bf Algorithm 4:} The RK-RRK algorithm for solving (\ref{mp}) with $\mbf b\in\ran(\mbf A\mbf B)$\\ 
\hline \noalign{\smallskip}
\quad {\bf Input}: $\mbf A\in\mbbr^{m\times \ell}$, $\mbf B\in\mbbr^{\ell\times n}$,  $\mbf b\in\mbbr^m$, $\gamma>0$, and maximum number of iterations {\tt maxit}.
\\\noalign{\smallskip}
\quad {\bf Output}: an approximation of the solution of (\ref{mp}) with $\mbf b\in\ran(\mbf A\mbf B)$.
\\ \noalign{\smallskip}
\quad {\bf Initialize}: $\mbf y^{(0)}=\mbf 0$, $\mbf z^{(0)}\in\ran(\mbf B^\top)$, and $\mbf x^{(0)}=\nabla f^*(\mbf z^{(0)})$.\\ \noalign{\smallskip}
\quad {\bf for} $k=1,2,\ldots,$ {\tt maxit} {\bf do}\\ \noalign{\smallskip}
\quad \qquad  Pick $j_k\in[m]$ with probability ${\|\mbf A_{j_k,:}\|^2_2}/{\|\mbf A\|_\rmf^2}$\\  \noalign{\smallskip}
\quad \qquad  Set $\mbf y^{(k)}=\mbf y^{(k-1)}-\dsp\frac{\mbf A_{j_k,:}\mbf y^{(k-1)}-b_{j_k}}{\|\mbf A_{j_k,:}\|_2^2}(\mbf A_{j_k,:})^\top$\\  \noalign{\smallskip}
\quad \qquad  Pick $i_k\in[\ell]$ with probability ${\|\mbf B_{i_k,:}\|^2_2}/{\|\mbf B\|_\rmf^2}$\\  \noalign{\smallskip}
\quad \qquad Set $\mbf z^{(k)} = \mbf z^{(k-1)}-\gamma\dsp \frac{\mbf B_{i_k,:}\mbf x^{(k-1)}-y^{(k)}_{i_k}}{\|\mbf B_{i_k,:}\|^2_2}(\mbf B_{i_k,:})^\top$ \\  \noalign{\smallskip}
\quad \qquad Set $\mbf x^{(k)}=\nabla f^*(\mbf z^{(k)})$ \\  \noalign{\smallskip}
\quad {\bf end}\\
\bottomrule
\end{tabular*}
\end{center}

\begin{theorem}\label{main} Let $\mbf x_\star$ be the solution of the problem {\rm(\ref{mp})} with $\mbf b\in\ran(\mbf A\mbf B)$. Assume that $f$ is $\gamma$-strongly convex and strongly admissible for $(\mbf A, \mbf B, \mbf b)$. Let $\nu>0$ be the constant  from {\rm(\ref{nu})}.  For any $\delta>0$, the sequences $\{\mbf x^{(k)}\}$ and $\{\mbf z^{(k)}\}$ generated by Algorithm $4$ satisfy \begin{align*}\mbbe\bem D_{f,\mbf z^{(k)}}(\mbf x^{(k)},\mbf x_\star)\eem  & \leq \l(1+\frac{\delta}{\gamma}\r)^k\beta^kD_{f,\mbf z^{(0)}}(\mbf x^{(0)},\mbf x_\star)+ \frac{(\delta+\gamma)\gamma}{2\delta\|\mbf B\|_\rmf^2}\|\mbf A^\dag\mbf b\|_2^2\sum_{i=0}^{k-1}\alpha^{k-i}\l(1+\frac{\delta}{\gamma}\r)^i\beta^i,\end{align*} where $$\alpha=1-\frac{\sigma_{\min}^2(\mbf A)}{\|\mbf A\|_\rmf^2},\quad \beta=1-\frac{\gamma\nu}{2\|\mbf B\|_\rmf^2}.$$\end{theorem}
\proof
Introduce the auxiliary vector $$\wh{\mbf z}^{(k)} := \mbf z^{(k-1)}- w(\mbf B_{i_k,:})^\top \quad \mbox{with}\quad w:= \gamma\frac{\mbf B_{i_k,:}\mbf x^{(k-1)}-\mbf I_{i_k,:}\mbf A^\dag\mbf b}{\|\mbf B_{i_k,:}\|^2_2}.$$  We have $$\mbf z^{(k)}-\wh{\mbf z}^{(k)}=\gamma\frac{y^{(k)}_{i_k}-\mbf I_{i_k,:}\mbf A^\dag\mbf b}{\|\mbf B_{i_k,:}\|^2_2}(\mbf B_{i_k,:})^\top.$$ Then, \beq\label{ywhy}
 	\|\mbf z^{(k)}-\wh{\mbf z}^{(k)}\|_2^2  = \gamma^2\frac{|y^{(k)}_{i_k}-\mbf I_{i_k,:}\mbf A^\dag\mbf b|^2}{\|\mbf B_{i_k,:}\|^2_2}.  \eeq
 	
Let $\mbbe_{k-1}\bem\cdot\eem$ denote the conditional expectation conditioned on $\mbf y^{(k-1)}$, $\mbf z^{(k-1)}$, and $\mbf x^{(k-1)}$. Let $\mbbe_{k-1}^{\rm B}\bem\cdot\eem$ denote the conditional expectation conditioned on $\mbf y^{(k)}$, $\mbf z^{(k-1)}$, and $\mbf x^{(k-1)}$. Then, by the law of total expectation, we have $$\mbbe_{k-1}\bem\cdot\eem=\mbbe_{k-1}\bem\mbbe_{k-1}^{\rm B}\bem\cdot\eem\eem.$$ Taking conditional expectation for (\ref{ywhy}) conditioned on $\mbf y^{(k-1)}$, $\mbf z^{(k-1)}$, and $\mbf x^{(k-1)}$, we obtain \begin{align*}\mbbe_{k-1}\bem\|\mbf z^{(k)}-\wh{\mbf z}^{(k)}\|_2^2\eem &=\mbbe_{k-1}\bem\mbbe_{k-1}^{\rm B}\bem\|\mbf z^{(k)}-\wh{\mbf z}^{(k)}\|_2^2\eem\eem\\ & = \frac{\gamma^2}{\|\mbf B\|_\rmf^2}\mbbe_{k-1}\bem\|\mbf y^{(k)}-\mbf A^\dag\mbf b\|_2^2 \eem. \end{align*} Then, by the law of total expectation and the estimate (\ref{az}), we have
\beq\mbbe\bem\|\mbf z^{(k)}-\wh{\mbf z}^{(k)}\|_2^2\eem = \frac{\gamma^2}{\|\mbf B\|_\rmf^2}\mbbe\bem\|\mbf y^{(k)}-\mbf A^\dag\mbf b\|_2^2 \eem  \leq \frac{\gamma^2\alpha^k}{\|\mbf B\|_\rmf^2}\|\mbf A^\dag\mbf b\|_2^2.\label{ey}\eeq

By $\mbf A\mbf B\mbf x_\star=\mbf b$ and $\rank(\mbf A)=\ell$, we have $\mbf B\mbf x_\star=\mbf A^\dag\mbf b$, which implies \beq\label{w}w=\gamma\frac{\mbf B_{i_k,:}\mbf x^{(k-1)}-\mbf I_{i_k,:}\mbf A^\dag\mbf b}{\|\mbf B_{i_k,:}\|^2_2}=\gamma\frac{\mbf B_{i_k,:}(\mbf x^{(k-1)}-\mbf x_\star)}{\|\mbf B_{i_k,:}\|^2_2}.\eeq
Define \beq\label{xhat}\wh{\mbf x}^{(k)}:=\nabla f^*(\wh{\mbf z}^{(k)}).\eeq By (\ref{xgz}), we have \beq\label{yhat}\wh{\mbf z}^{(k)}\in\p f(\wh{\mbf x}^{(k)}).\eeq 
The Bregman distance between $\wh{\mbf x}^{(k)}$ and $\mbf x_\star$ with respect to $f$ and $\wh{\mbf z}^{(k)}$ satisfies \begin{align*}
D_{f,\wh{\mbf z}^{(k)}}(\wh{\mbf x}^{(k)},\mbf x_\star) &\stackrel{(\ref{dfz})}= f^*(\wh{\mbf z}^{(k)})-\la\wh{\mbf z}^{(k)},\mbf x_\star\ra+f(\mbf x_\star) \\ & = f^*( \mbf z^{(k-1)}-w(\mbf B_{i_k,:})^\top) -\la \mbf z^{(k-1)}-w(\mbf B_{i_k,:})^\top,\mbf x_\star\ra+f(\mbf x_\star) \\ & \stackrel{(\ref{sc2})}\leq f^*(\mbf z^{(k-1)}) -\la\nabla f^*(\mbf z^{(k-1)}),w(\mbf B_{i_k,:})^\top \ra+\frac{1}{2\gamma}\|w(\mbf B_{i_k,:})^\top\|_2^2\\ &\qquad	-\la \mbf z^{(k-1)}-w(\mbf B_{i_k,:})^\top,\mbf x_\star\ra+f(\mbf x_\star) \\ & \stackrel{(\ref{dfz})}= D_{f,\mbf z^{(k-1)}}(\mbf x^{(k-1)},\mbf x_\star)-\la\mbf x^{(k-1)},w(\mbf B_{i_k,:})^\top\ra\\ &\qquad +\frac{1}{2\gamma}\|w(\mbf B_{i_k,:})^\top\|_2^2+\la(w\mbf B_{i_k,:})^\top,\mbf x_\star\ra\\ & =  D_{f,\mbf z^{(k-1)}}(\mbf x^{(k-1)},\mbf x_\star)-\la\mbf x^{(k-1)}-\mbf x_\star,w(\mbf B_{i_k,:})^\top\ra+\frac{1}{2\gamma}\|w(\mbf B_{i_k,:})^\top\|_2^2\\ & \stackrel{(\ref{w})}  = D_{f,\mbf z^{(k-1)}}(\mbf x^{(k-1)},\mbf x_\star)-\gamma\frac{|\mbf B_{i_k,:}(\mbf x^{(k-1)}-\mbf x_\star)|_2^2}{\|\mbf B_{i_k,:}\|_2^2} +\frac{\gamma^2}{2\gamma}\frac{|\mbf B_{i_k,:}(\mbf x^{(k-1)}-\mbf x_\star)|_2^2}{\|\mbf B_{i_k,:}\|_2^2} \\ &= D_{f,\mbf z^{(k-1)}}(\mbf x^{(k-1)},\mbf x_\star)-\frac{\gamma}{2}\frac{|\mbf B_{i_k,:}(\mbf x^{(k-1)}-\mbf x_\star)|_2^2}{\|\mbf B_{i_k,:}\|_2^2}.
\end{align*} By induction, we can prove that $\mbf z^{(k)}\in\ran(\mbf B^\top)$. By $\mbf x^{(k)}=\nabla f^*(\mbf z^{(k)})$ and (\ref{xgz}), we have $\mbf z^{(k)}\in\p f(\mbf x^{(k)})$. 
Taking conditional expectation conditioned on $\mbf y^{(k-1)}$, $\mbf z^{(k-1)}$, and $\mbf x^{(k-1)}$, we have 
\begin{align*}
\mbbe_{k-1}\bem D_{f,\wh{\mbf z}^{(k)}}(\wh{\mbf x}^{(k)},\mbf x_\star) \eem & \leq D_{f,\mbf z^{(k-1)}}(\mbf x^{(k-1)},\mbf x_\star) -\frac{\gamma}{2\|\mbf B\|_\rmf^2}\|\mbf B(\mbf x^{(k-1)}-\mbf x_\star)\|_2^2\\ &\stackrel{(\ref{nu})}\leq \l(1-\frac{\gamma\nu}{2\|\mbf B\|_\rmf^2}\r)D_{f,\mbf z^{(k-1)}}(\mbf x^{(k-1)},\mbf x_\star).
\end{align*}
Thus, by the law of total expectation, we have
\beq\label{ed} \mbbe \bem D_{f,\wh{\mbf z}^{(k)}}(\wh{\mbf x}^{(k)},\mbf x_\star)\eem\leq \beta\mbbe\bem D_{f,\mbf z^{(k-1)}}(\mbf x^{(k-1)},\mbf x_\star)\eem.\eeq

Now, we consider the Bregman distance $D_{f,\mbf z^{(k)}}(\mbf x^{(k)},\mbf x_\star)$, which satisfies \begin{align*} D_{f,\mbf z^{(k)}}(\mbf x^{(k)},\mbf x_\star) &\stackrel{(\ref{dfz})} =D_{f,\wh{\mbf z}^{(k)}}(\wh{\mbf x}^{(k)},\mbf x_\star)+f^*(\mbf z^{(k)})-f^*(\wh{\mbf z}^{(k)})-\la\mbf z^{(k)},\mbf x_\star\ra+\la\wh{\mbf z}^{(k)},\mbf x_\star\ra\\
& \stackrel{(\ref{sc2})}\leq D_{f,\wh{\mbf z}^{(k)}}(\wh{\mbf x}^{(k)},\mbf x_\star)+\la\nabla f^*(\wh{\mbf z}^{(k)}),\mbf z^{(k)}-\wh{\mbf z}^{(k)}\ra+\frac{1}{2\gamma}\|\mbf z^{(k)}-\wh{\mbf z}^{(k)}\|_2^2\\ &\quad -\la\mbf z^{(k)}-\wh{\mbf z}^{(k)},\mbf x_\star\ra
\\ & \stackrel{(\ref{xhat})}= D_{f,\wh{\mbf z}^{(k)}}(\wh{\mbf x}^{(k)},\mbf x_\star)+\la\wh{\mbf x}^{(k)}-\mbf x_\star,\mbf z^{(k)}-\wh{\mbf z}^{(k)}\ra+\frac{1}{2\gamma}\|\mbf z^{(k)}-\wh{\mbf z}^{(k)}\|_2^2
\\ & \leq  D_{f,\wh{\mbf z}^{(k)}}(\wh{\mbf x}^{(k)},\mbf x_\star)+\frac{\delta}{2}\|\wh{\mbf x}^{(k)}-\mbf x_\star\|_2^2+\frac{1}{2\delta}\|\mbf z^{(k)}-\wh{\mbf z}^{(k)}\|_2^2\\ &\qquad +\frac{1}{2\gamma}\|\mbf z^{(k)}-\wh{\mbf z}^{(k)}\|_2^2 \quad \mbox{(by  Cauchy--Schwarz and Young's inequality) }
\\ &\stackrel{(\ref{yhat})(\ref{gamma})}\leq \l(1+\frac{\delta}{\gamma}\r)D_{f,\wh{\mbf z}^{(k)}}(\wh{\mbf x}^{(k)},\mbf x_\star)+\frac{\delta+\gamma}{2\delta\gamma}\|\mbf z^{(k)}-\wh{\mbf z}^{(k)}\|_2^2.
\end{align*}
Taking expectation, we have
\begin{align*}
	\mbbe\bem D_{f,\mbf z^{(k)}}(\mbf x^{(k)},\mbf x_\star)\eem & \leq  \l(1+\frac{\delta}{\gamma}\r)\mbbe\bem D_{f,\wh{\mbf z}^{(k)}}(\wh{\mbf x}^{(k)},\mbf x_\star)\eem+\frac{\delta+\gamma}{2\delta\gamma}\mbbe\bem\|\mbf z^{(k)}-\wh{\mbf z}^{(k)}\|_2^2\eem
	\\ & \stackrel{(\ref{ey})(\ref{ed})} \leq \l(1+\frac{\delta}{\gamma}\r)\beta\mbbe\bem D_{f,\mbf z^{(k-1)}}(\mbf x^{(k-1)},\mbf x_\star)\eem+\frac{(\delta+\gamma)\gamma}{2\delta}\frac{\alpha^k}{\|\mbf B\|_\rmf^2}\|\mbf A^\dag\mbf b\|_2^2
	\\ & \leq \l(1+\frac{\delta}{\gamma}\r)^2\beta^2\mbbe\bem D_{f,\mbf z^{(k-2)}}(\mbf x^{(k-2)},\mbf x_\star)\eem\\ &\quad  +\frac{(\delta+\gamma)\gamma}{2\delta\|\mbf B\|_\rmf^2}\|\mbf A^\dag\mbf b\|_2^2\l(\alpha^k+\alpha^{k-1}\l(1+\frac{\delta}{\gamma}\r)\beta\r)
	\\ & \leq \cdots 
	\\ & \leq \l(1+\frac{\delta}{\gamma}\r)^k\beta^kD_{f,\mbf z^{(0)}}(\mbf x^{(0)},\mbf x_\star)+\frac{(\delta+\gamma)\gamma}{2\delta\|\mbf B\|_\rmf^2}\|\mbf A^\dag\mbf b\|_2^2\sum_{i=0}^{k-1}\alpha^{k-i}\l(1+\frac{\delta}{\gamma}\r)^i\beta^i.
\end{align*} This completes the proof.
\endproof

\begin{remark}\label{remark}
Let $\rho=\max\{\alpha,\beta\}$. For all $\delta$ satisfying $(1+\delta/\gamma)\rho<1$, from Theorem {\rm \ref{main}}, we have \begin{align*}\mbbe\bem D_{f,\mbf z^{(k)}}(\mbf x^{(k)},\mbf x_\star)\eem 
& \leq \l(1+\frac{\delta}{\gamma}\r)^k\rho^kD_{f,\mbf z^{(0)}}(\mbf x^{(0)},\mbf x_\star)+\frac{(\delta+\gamma)\gamma}{2\delta\|\mbf B\|_\rmf^2}\|\mbf A^\dag\mbf b\|_2^2\rho^k\sum_{i=0}^{k-1}\l(1+\frac{\delta}{\gamma}\r)^i \\ 
& \leq \l(1+\frac{\delta}{\gamma}\r)^k\rho^kD_{f,\mbf z^{(0)}}(\mbf x^{(0)},\mbf x_\star)+\frac{(\delta+\gamma)\gamma}{2\delta\|\mbf B\|_\rmf^2}\|\mbf A^\dag\mbf b\|_2^2\rho^k\frac{\gamma}{\delta}\l(\l(1+\frac{\delta}{\gamma}\r)^k-1\r) \\
&\leq \l(1+\frac{\delta}{\gamma}\r)^k\rho^k\l(D_{f,\mbf z^{(0)}}(\mbf x^{(0)},\mbf x_\star)+\frac{(\delta+\gamma)\gamma^2}{2\delta^2\|\mbf B\|_\rmf^2}\|\mbf A^\dag\mbf b\|_2^2\r).\end{align*} Using \eqref{gamma}, we have 
$$\mbbe\bem\|\mbf x^{(k)}-\mbf x_\star\|_2^2\eem\leq \l(1+\frac{\delta}{\gamma}\r)^k\rho^k\frac{2}{\gamma}\l(D_{f,\mbf z^{(0)}}(\mbf x^{(0)},\mbf x_\star)+\frac{(\delta+\gamma)\gamma^2}{2\delta^2\|\mbf B\|_\rmf^2}\|\mbf A^\dag\mbf b\|_2^2\r),$$ 
which means that the RK-RRK algorithm converges linearly in expectation to the solution of the problem \eqref{mp} with the rate $(1+\delta/\gamma)\rho$.	
\end{remark}

In Algorithm 4 we use $\mbf y^{(0)}=\mbf 0$ for simplicity, and the analysis for any $\mbf y^{(0)}\in\mbbr^\ell$ is straightforward. For the choice $\dsp f(\mbf x)=\frac{1}{2}\|\mbf x\|_2^2$, Algorithm 4 becomes the RK-RK algorithm \cite{ma2018itera}. For the choice $\dsp f(\mbf x)=\frac{1}{2}\|\mbf x\|_2^2+\lambda\|\mbf x\|_1$, we call the resulting algorithm the RK-RSK algorithm. In the following remark, we give the relationship between the GERK-(a,d) algorithm \cite{schopfer2022exten} and the RK-RSK algorithm. Recall that $\dsp f(\mbf x)=\frac{1}{2}\|\mbf x\|_2^2+\lambda\|\mbf x\|_1$ is $1$-strongly convex.

\begin{remark} The iterates of the {\rm GERK-(a,d)} algorithm {\rm\cite{schopfer2022exten}} for sparse (least squares) solutions of the full linear system $\mbf C\mbf x=\mbf b$ are \begin{align*}\mbf y^{(k)}&=\mbf y^{(k-1)}-\dsp\frac{(\mbf C_{:,j_k})^\top\mbf y^{(k-1)}}{\|\mbf C_{:,j_k}\|_2^2}\mbf C_{:,j_k},\\ \mbf z^{(k)} &= \mbf z^{(k-1)}-\dsp \frac{\mbf C_{i_k,:}\mbf x^{(k-1)}-b_{i_k}+y^{(k)}_{i_k}}{\|\mbf C_{i_k,:}\|^2_2}(\mbf C_{i_k,:})^\top, \\ \mbf x^{(k)}&=S_\lambda(\mbf z^{(k)}),\end{align*} with initial iterates $\mbf y^{(0)}=\mbf b$, $\mbf z^{(0)}\in\ran(\mbf C^\top)$, and $\mbf x^{(0)}=S_\lambda(\mbf z^{(0)})$. Note that the normal equations $\mbf C^\top\mbf C\mbf x=\mbf C^\top\mbf b$ can be viewed as the factorized linear system $\wh{\mbf A}\wh{\mbf B}\mbf x=\wh{\mbf b}$ with $\wh{\mbf A}=\mbf C^\top$, $\wh{\mbf B}=\mbf C$, and $\wh{\mbf b}=\mbf C^\top\mbf b$. The iterates of the {\rm RK-RSK} algorithm for $\wh{\mbf A}\wh{\mbf B}\mbf x=\wh{\mbf b}$ are \begin{align*}\wh{\mbf y}^{(k)}&=\wh{\mbf y}^{(k-1)}-\dsp\frac{\wh{\mbf A}_{j_k,:}\wh{\mbf y}^{(k-1)}-\wh{b}_{j_k}}{\|\wh{\mbf A}_{j_k,:}\|_2^2}(\wh{\mbf A}_{j_k,:})^\top,\\ \wh{\mbf z}^{(k)} &= \wh{\mbf z}^{(k-1)}-\dsp \frac{\wh{\mbf B}_{i_k,:}\wh{\mbf x}^{(k-1)}-\wh{y}^{(k)}_{i_k}}{\|\wh{\mbf B}_{i_k,:}\|^2_2}(\wh{\mbf B}_{i_k,:})^\top, \\ \wh{\mbf x}^{(k)}&=S_\lambda(\wh{\mbf z}^{(k)}),\end{align*} with initial iterates $\wh{\mbf y}^{(0)}=\mbf 0$, $\mbf z^{(0)}\in\ran(\wh{\mbf B}^\top)$, and $\wh{\mbf x}^{(0)}=S_\lambda(\wh{\mbf z}^{(0)})$. By using $\wh{\mbf A}_{j_k,:}=(\mbf C_{:,j_k})^\top$, $\wh {\mbf B}_{i_k,:}=\mbf C_{i_k,:}$, and $\wh{b}_{j_k}=(\mbf C_{:,j_k})^\top\mbf b$, we have \begin{align*}\wh{\mbf y}^{(k)}&=\wh{\mbf y}^{(k-1)}-\dsp\frac{(\mbf C_{:,j_k})^\top\wh{\mbf y}^{(k-1)}-(\mbf C_{:,j_k})^\top\mbf b}{\|\mbf C_{:,j_k}\|_2^2}\mbf C_{:,j_k},\\ \wh{\mbf z}^{(k)} &= \wh{\mbf z}^{(k-1)}-\dsp \frac{\mbf C_{i_k,:}\wh{\mbf x}^{(k-1)}-\wh{y}^{(k)}_{i_k}}{\|\wh{\mbf B}_{i_k,:}\|^2_2}(\mbf C_{i_k,:})^\top.\end{align*} It follows that $$\mbf b-\wh{\mbf y}^{(k)}=\mbf b-\wh{\mbf y}^{(k-1)}-\dsp\frac{(\mbf C_{:,j_k})^\top(\mbf b-\wh{\mbf y}^{(k-1)})}{\|\mbf C_{:,j_k}\|_2^2}\mbf C_{:,j_k}.$$ If $\mbf z^{(0)}=\wh{\mbf z}^{(0)}$, then by induction it is straightforward to prove that the iterates $\mbf x^{(k)}$, $\mbf y^{(k)}$, and $\mbf z^{(k)}$ of the {\rm GERK-(a,d)} algorithm for $\mbf C\mbf x=\mbf b$ are equal to $\wh{\mbf x}^{(k)}$, $\mbf b-\wh{\mbf y}^{(k)}$, and $\wh{\mbf z}^{(k)}$, respectively. 
\end{remark}

\subsection{The RGS-RRK algorithm for the case $\mbf b\notin\ran(\mbf A\mbf B)$}

In this subsection, we analyze the RGS-RRK algorithm (Algorithm 5) and prove its linear convergence property. We note that $\mbf y^{(k)}$ in the RGS-RRK algorithm is actually the $k$th iterate of the RGS algorithm for $\min_{\mbf y}\|\mbf b-\mbf A\mbf y\|_2$, whose convergence estimate is given in (\ref{rgsbound}). We mention that the auxiliary vector $\mbf r^{(k)} = \mbf b-\mbf A\mbf y^{(k)}$ with the update rule $\mbf r^{(k)}=\mbf r^{(k-1)}-d_k\mbf A_{:,j_k}$ is introduced to avoid the computation of the matrix-vector multiplication $\mbf A\mbf y^{(k)}$. The vectors $\mbf z^{(k)}$ and $\mbf x^{(k)}$ are one-step updates of the RRK algorithm for $\mbf B\mbf x=\mbf y^{(k)}$ from $\mbf z^{(k-1)}$ and $\mbf x^{(k-1)}$. We give the convergence result of the RGS-RRK algorithm in Theorem \ref{main2}. 

\begin{center}
\begin{tabular*}{160mm}{l}
\toprule {\bf Algorithm 5:} The RGS-RRK algorithm for solving (\ref{mp}) with $\mbf b\notin\ran(\mbf A\mbf B)$\\ 
\hline \noalign{\smallskip}
\quad {\bf Input}: $\mbf A\in\mbbr^{m\times \ell}$, $\mbf B\in\mbbr^{\ell\times n}$,  $\mbf b\in\mbbr^m$, $\gamma>0$, and maximum number of  iterations {\tt maxit}.
\\ \noalign{\smallskip}
\quad {\bf Output}: an approximation of the solution of (\ref{mp}) with $\mbf b\notin\ran(\mbf A\mbf B)$.
\\ \noalign{\smallskip}
\quad {\bf Initialize}: $\mbf y^{(0)}=\mbf 0$, $\mbf r^{(0)}=\mbf b$, $\mbf z^{(0)}\in\ran(\mbf B^\top)$, and $\mbf x^{(0)}=\nabla f^*(\mbf z^{(0)})$.\\ \noalign{\smallskip}
\quad {\bf for} $k=1,2,\ldots,$ {\tt maxit} {\bf do}\\ \noalign{\smallskip}
\quad \qquad Pick $j_k\in[\ell]$ with probability ${\|\mbf A_{:,j_k}\|^2_2}/{\|\mbf A\|_\rmf^2}$\\  \noalign{\smallskip}
\quad \qquad Compute $\dsp d_k=\frac{(\mbf A_{:,j_k})^\top\mbf r^{(k-1)}}{\|\mbf A_{:,j_k}\|_2^2}$\\  \noalign{\smallskip}
\quad \qquad  Set $y_{j_k}^{(k)}= y_{j_k}^{(k-1)}+d_k$, $y_j^{(k)}= y_j^{(k-1)}$ for $j\neq j_k$, and $\mbf r^{(k)}=\mbf r^{(k-1)}-d_k\mbf A_{:,j_k}$\\  \noalign{\smallskip}
\quad \qquad  Pick $i_k\in[\ell]$ with probability ${\|\mbf B_{i_k,:}\|^2_2}/{\|\mbf B\|_\rmf^2}$\\  \noalign{\smallskip}
\quad \qquad Set $\mbf z^{(k)} = \mbf z^{(k-1)}-\gamma\dsp \frac{\mbf B_{i_k,:}\mbf x^{(k-1)}-y^{(k)}_{i_k}}{\|\mbf B_{i_k,:}\|^2_2}(\mbf B_{i_k,:})^\top$ \\  \noalign{\smallskip}
\quad \qquad Set $\mbf x^{(k)}=\nabla f^*(\mbf z^{(k)})$ \\  \noalign{\smallskip}
\quad {\bf end}\\
\bottomrule
\end{tabular*}
\end{center}

\begin{theorem}\label{main2} Let $\mbf x_\star$ be the solution of the problem {\rm(\ref{mp})} with $\mbf b\notin\ran(\mbf A\mbf B)$. Assume that $f$ is $\gamma$-strongly convex and strongly admissible for $(\mbf A, \mbf B, \mbf b)$. Let $\nu>0$ be the constant  from {\rm(\ref{nu})}.  For any $\delta>0$, the sequences $\{\mbf x^{(k)}\}$ and $\{\mbf z^{(k)}\}$ generated by Algorithm $5$ satisfy \begin{align*}\mbbe\bem D_{f,\mbf z^{(k)}}(\mbf x^{(k)},\mbf x_\star)\eem  & \leq \frac{(\delta+\gamma)\gamma}{2\delta\|\mbf B\|_\rmf^2}\|\mbf A^\dag\|_2^2\|\mbf A\mbf A^\dag\mbf b\|_2^2\sum_{i=0}^{k-1}\alpha^{k-i}\l(1+\frac{\delta}{\gamma}\r)^i\beta^i \\ &\quad  +\l(1+\frac{\delta}{\gamma}\r)^k\beta^kD_{f,\mbf z^{(0)}}(\mbf x^{(0)},\mbf x_\star),\end{align*} where $$\alpha=1-\frac{\sigma_{\min}^2(\mbf A)}{\|\mbf A\|_\rmf^2},\quad \beta=1-\frac{\gamma\nu}{2\|\mbf B\|_\rmf^2}.$$\end{theorem}
\proof
The proof of Theorem \ref{main2} is almost the same as that of Theorem \ref{main} except that the estimate (\ref{ey}) is replaced with $$\mbbe\bem\|\mbf z^{(k)}-\wh{\mbf z}^{(k)}\|_2^2\eem = \frac{\gamma^2}{\|\mbf B\|_\rmf^2}\mbbe\bem\|\mbf y^{(k)}-\mbf A^\dag\mbf b\|_2^2 \eem \stackrel{(\ref{rgsbound})}\leq \frac{\gamma^2\alpha^k}{\|\mbf B\|_\rmf^2}\|\mbf A^\dag\|_2^2\|\mbf A\mbf A^\dag\mbf b\|_2^2.$$ We omit the  details. 
\endproof

By the similar discussion as in Remark \ref{remark}, we obtain that the RGS-RRK algorithm converges linearly in expectation to the solution of the problem (\ref{mp}). In Algorithm 5 we use $\mbf y^{(0)}=\mbf 0$ for simplicity, and the analysis for any $\mbf y^{(0)}\in\mbbr^\ell$ is straightforward. For the choice $\dsp f(\mbf x)=\frac{1}{2}\|\mbf x\|_2^2$, Algorithm 5 becomes the RGS-RK algorithm \cite{zhao2023rando}. For the choice $\dsp f(\mbf x)=\frac{1}{2}\|\mbf x\|_2^2+\lambda\|\mbf x\|_1$, we call the resulting algorithm the RGS-RSK algorithm.

\begin{remark}
	The idea of the {\rm RGS-RSK} algorithm can also be used to design a randomized sparse extended Gauss--Seidel {\rm(RSEGS)} algorithm for sparse $($least squares$)$ solutions of the full linear system $\mbf C\mbf x=\mbf b$. The iterates of the proposed {\rm RSEGS} algorithm are \begin{align*}\mbf y^{(k)}&=\mbf y^{(k-1)}-\dsp\frac{(\mbf C_{:,j_k})^\top(\mbf C\mbf y^{(k-1)}-\mbf b)}{\|\mbf C_{:,j_k}\|_2^2}\mbf I_{:,j_k},\\ \mbf z^{(k)} &= \mbf z^{(k-1)}-\dsp \frac{\mbf C_{i_k,:}(\mbf x^{(k-1)}-\mbf y^{(k)})}{\|\mbf C_{i_k,:}\|^2_2}(\mbf C_{i_k,:})^\top, \\ \mbf x^{(k)}&=S_\lambda(\mbf z^{(k)}),\end{align*} with initial iterates $\mbf y^{(0)}\in\mbbr^n$, $\mbf z^{(0)}\in\ran(\mbf C^\top)$, and $\mbf x^{(0)}=\nabla f^*(\mbf z^{(0)})$.  Here, $\mbf y^{(k)}$ is the $k$th iterate of the {\rm RGS} algorithm for $\min_{\mbf y}\|\mbf b-\mbf C\mbf y\|_2$, and $\mbf z^{(k)}$ and $\mbf x^{(k)}$ are one-step updates of the {\rm RSK} algorithm for $\mbf C \mbf x=\mbf C\mbf y^{(k)}$ from $\mbf z^{(k-1)}$ and $\mbf x^{(k-1)}$. By using the same technique used in the proof of Theorem $\ref{main2}$, we can show that the {\rm RSEGS} algorithm linearly converges to the unique solution of the minimization problem  \eqref{cls}.
\end{remark}

\section{Computed examples}

In this section, we report some numerical results for the proposed algorithms for sparse (least squares) solutions of factorized linear systems. All experiments are performed using MATLAB R2020b on a laptop with 2.7 GHz Quad-Core Intel Core i7 processor, 16 GB memory, and Mac operating system.   
 
\subsection{Example 1}
The matrices $\mbf A$ and $\mbf B$ are Gaussian matrices generated by {\tt A=randn($m,\ell$)}  and  {\tt B=randn($\ell,n$)}. We construct a sparse vector $\mbf x_\star$ with $s$ normally distributed  non-zero entries, whose support is randomly generated. Then we set $\mbf b=\mbf A\mbf B\mbf x_\star$ for $\mbf b\in\ran(\mbf A\mbf B)$, and set $\mbf b=\wh{\mbf b}+\wh{\mbf b}_\bot$ for $\mbf b\notin\ran(\mbf A\mbf B)$ with  $\wh{\mbf b}=\mbf A\mbf B\mbf x_\star$ and $\wh{\mbf b}_\bot=\mbf N\mbf v\|\wh{\mbf b}\|_2/\|\mbf N\mbf v\|_2\in\nul(\mbf B^\top\mbf A^\top)=\nul(\mbf A^\top)$, where the columns of $\mbf N$ form an orthonormal basis of $\nul(\mbf A^\top)$ and $\mbf v$ is a Gaussian vector generated by {\tt v=randn($m-\ell$,1)}. For the case $\mbf b\in\ran(\mbf A\mbf B)$, we compare the proposed RK-RSK algorithm with the RK-RK algorithm \cite{ma2018itera}. For the case $\mbf b\notin\ran(\mbf A\mbf B)$, we compare the proposed RGS-RSK algorithm with the RGS-RK algorithm \cite{zhao2023rando}. For the proposed algorithms, we use $\lambda=1$, $\mbf y^{(0)}=\mbf 0$, $\mbf z^{(0)} = \mbf 0$, and the maximum number of iterations {\tt maxit=20$m$}.

\begin{figure}[htb]
\centerline{\epsfig{figure=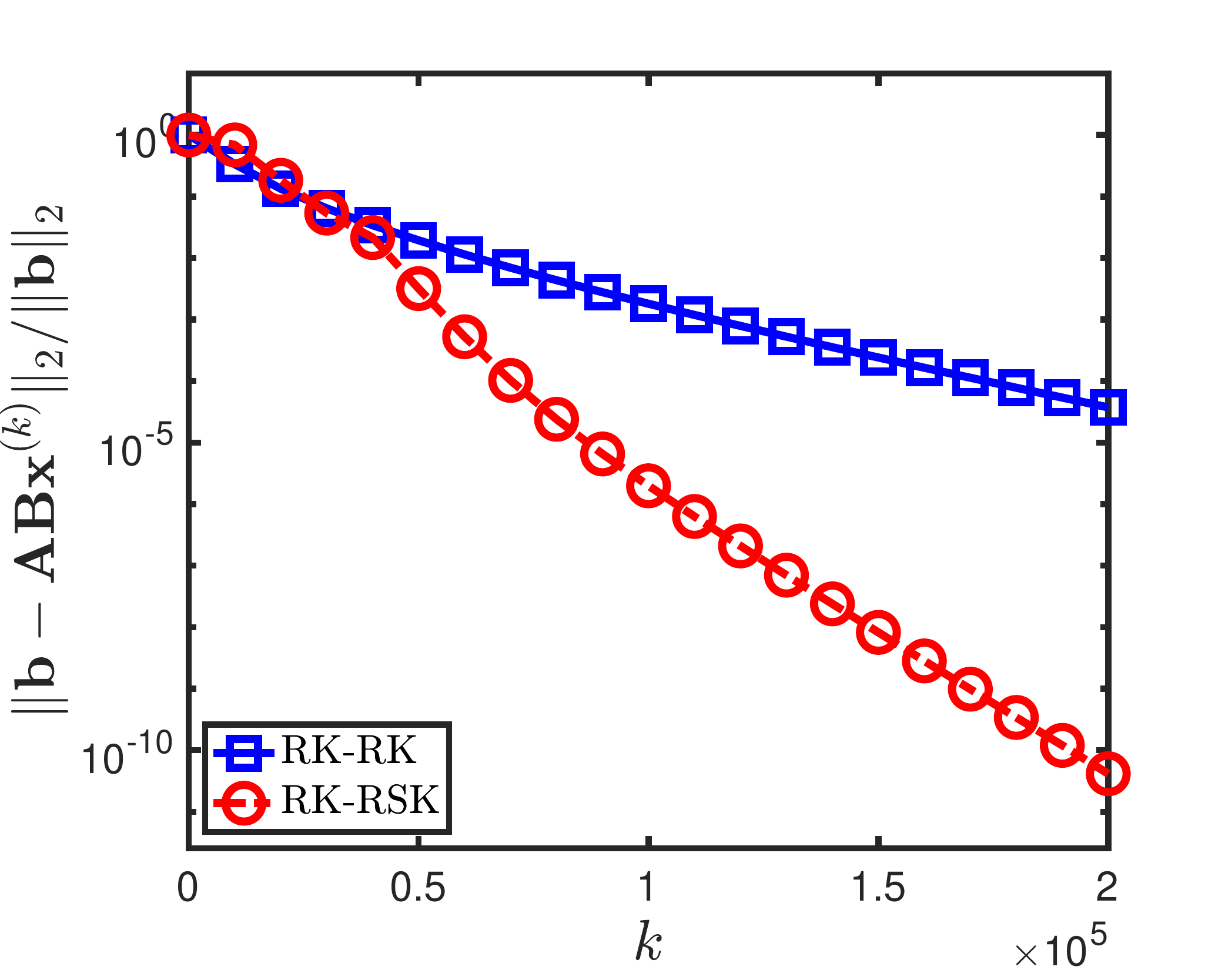,height=1.7in}\epsfig{figure=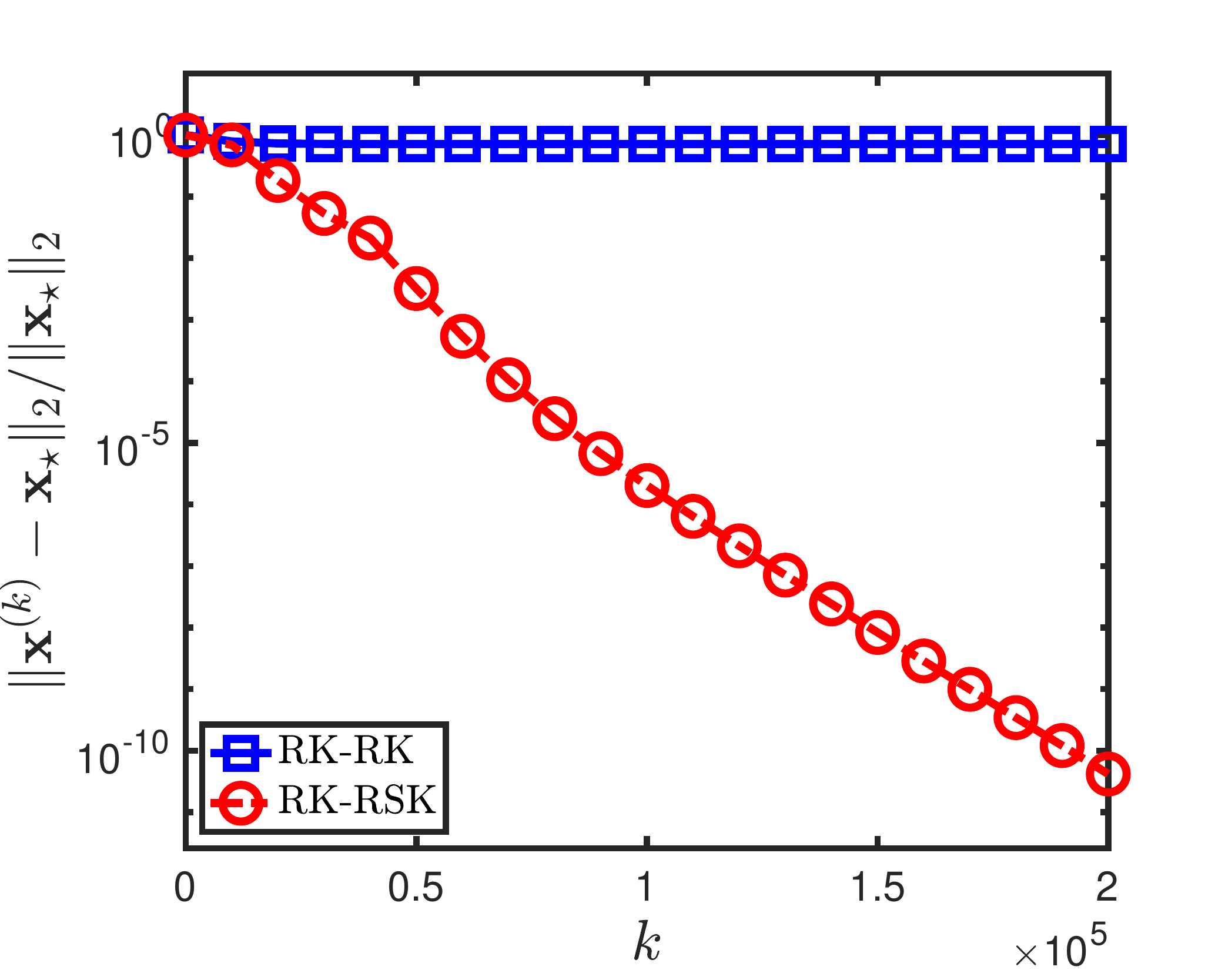,height=1.7in}\epsfig{figure=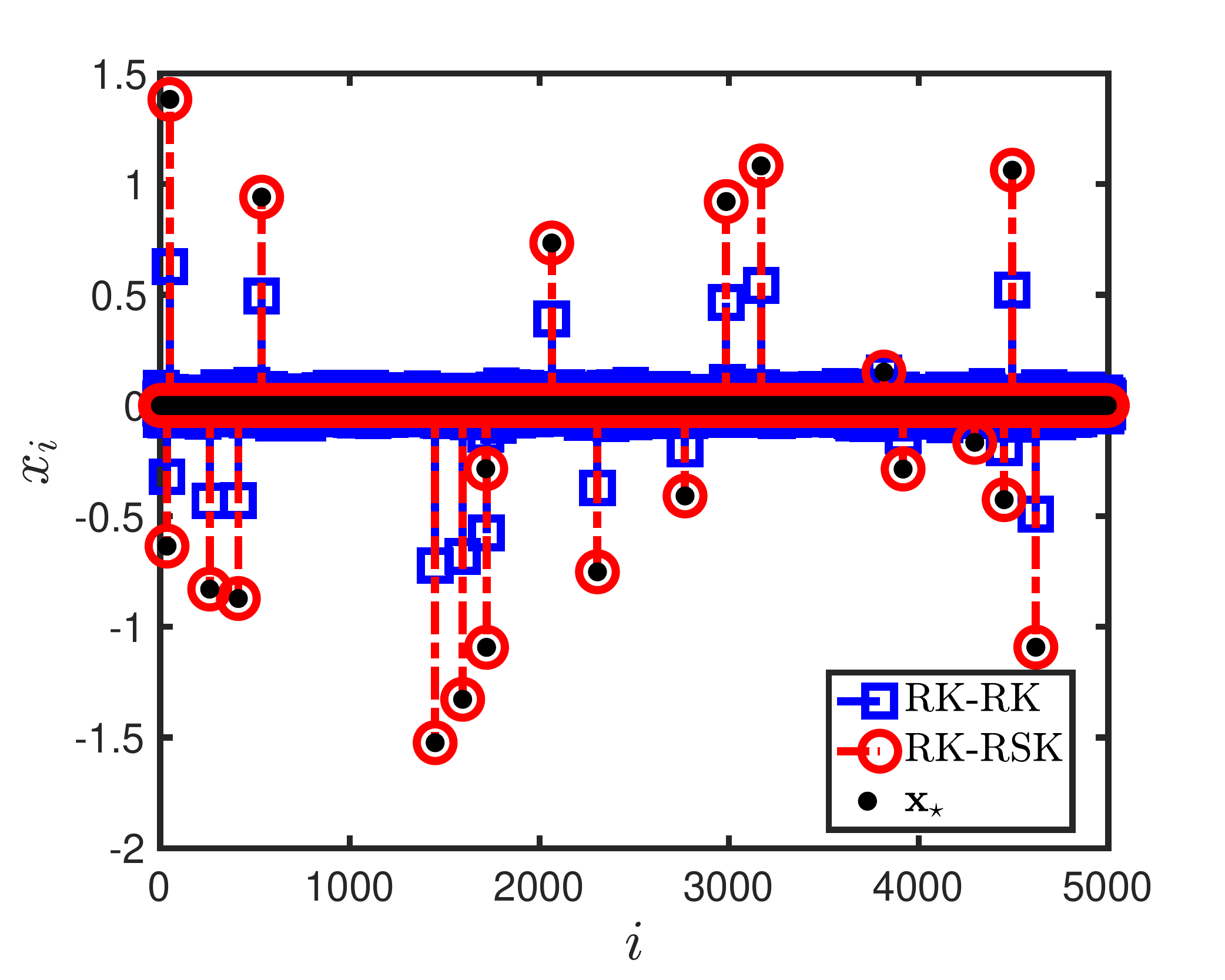,height=1.7in}}
\caption{Comparison of RK-RK and RK-RSK for the case $\mbf b\in\ran(\mbf A\mbf B)$. Left: iteration vs relative residual $\|\mbf b-\mbf A\mbf B\mbf x^{(k)}\|_2/\|\mbf b\|_2$. Middle: iteration vs relative error $\|\mbf x^{(k)}-\mbf x_\star\|_2/\|\mbf x_\star\|_2$. Right: approximated solutions (last iterates of RK-RK and RK-RSK) and the sparse vector $\mbf x_\star$.}
\label{fig1}
\end{figure} 

\begin{figure}[htb]
\centerline{\epsfig{figure=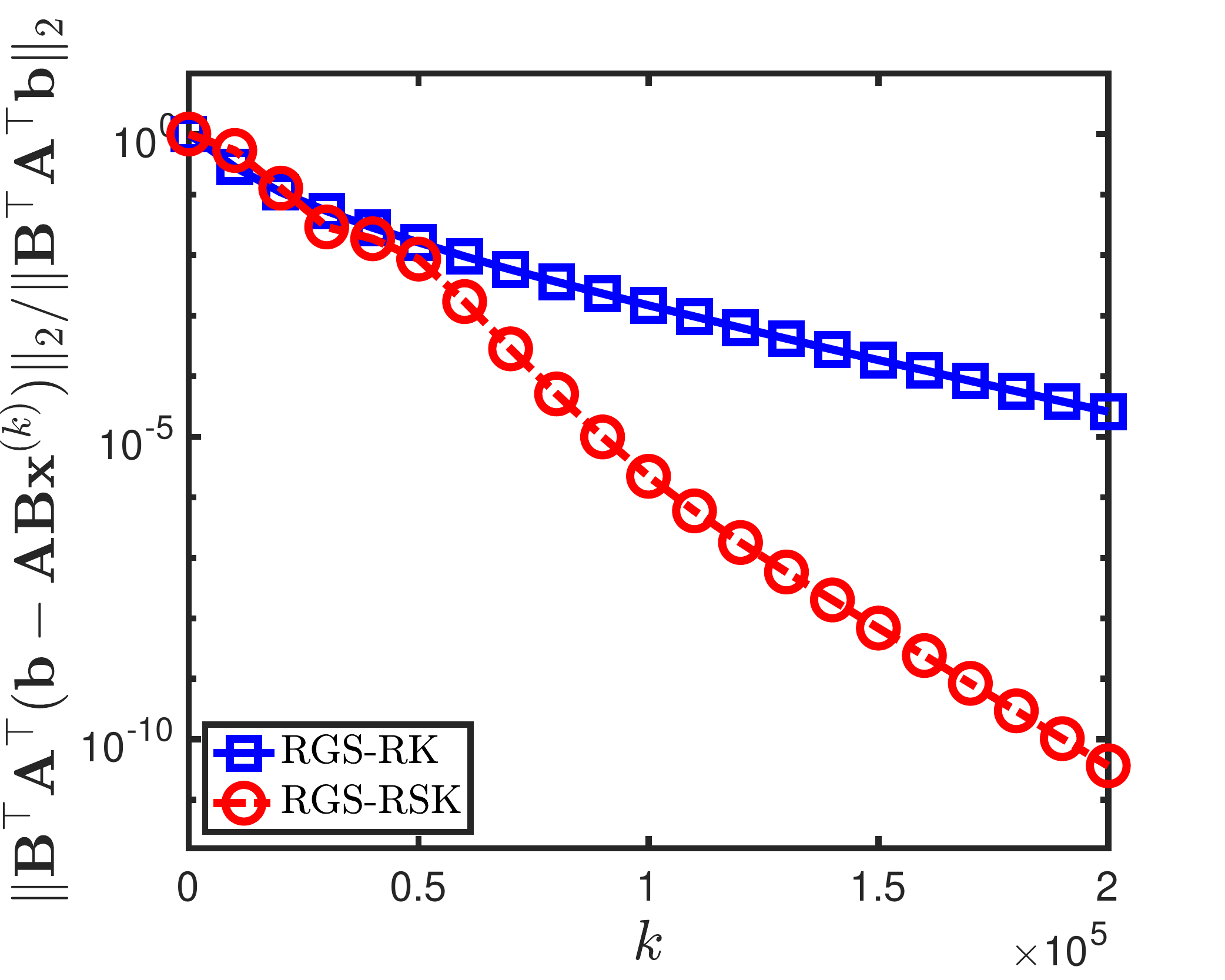,height=1.7in}\epsfig{figure=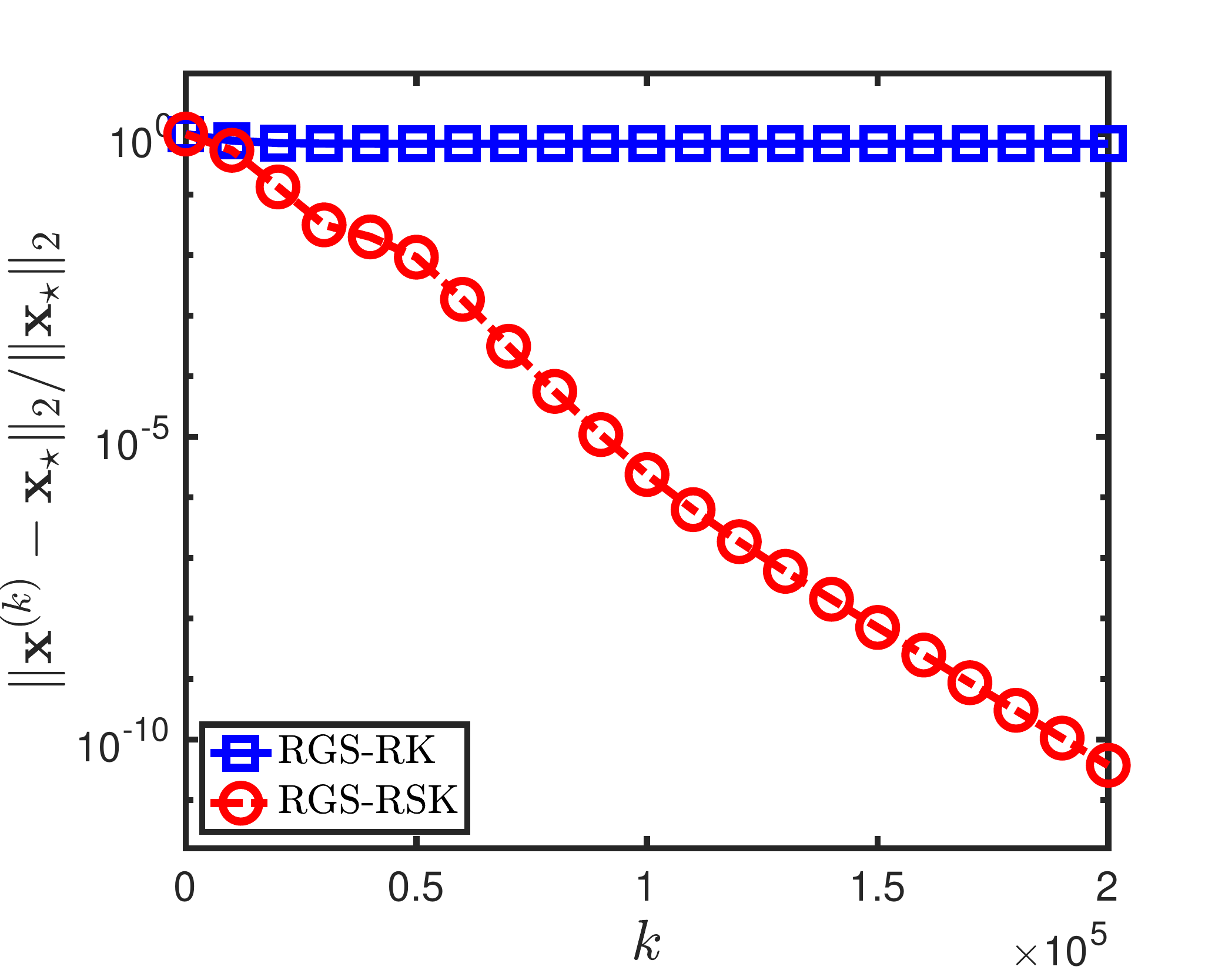,height=1.7in}\epsfig{figure=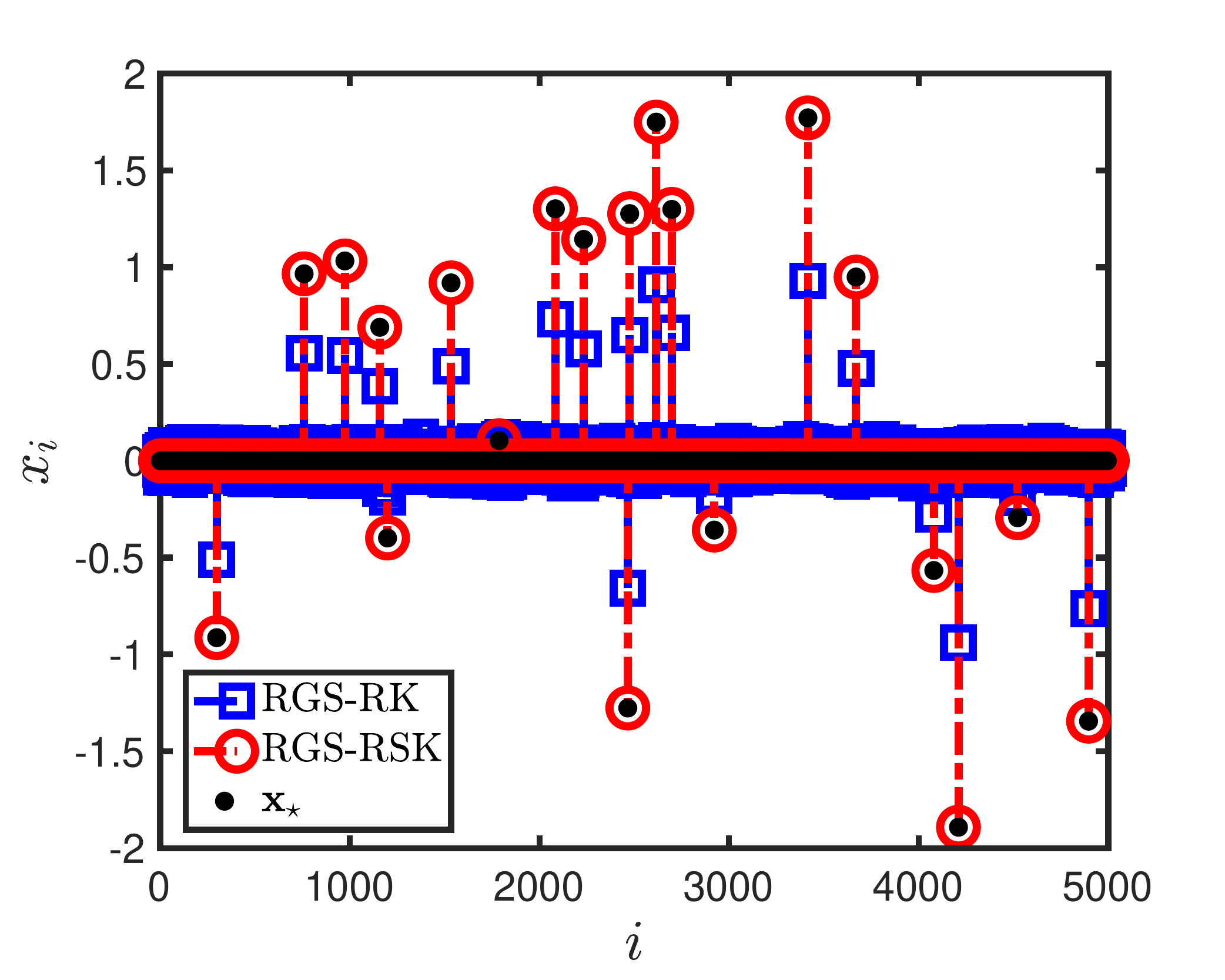,height=1.7in}}
\caption{Comparison of RGS-RK and RGS-RSK for the case $\mbf b\notin\ran(\mbf A\mbf B)$. Left: iteration vs relative residual $\|\mbf B^\top\mbf A^\top(\mbf b-\mbf A\mbf B\mbf x^{(k)})\|_2/\|\mbf B^\top\mbf A^\top\mbf b\|_2$. Middle: iteration vs relative error $\|\mbf x^{(k)}-\mbf x_\star\|_2/\|\mbf x_\star\|_2$. Right: approximated solutions (last iterates of RGS-RK and RGS-RSK) and the sparse vector $\mbf x_\star$.}
\label{fig2}
\end{figure} 

 In Figure \ref{fig1}, we plot the results for the case $\mbf b\in\ran(\mbf A\mbf B)$ with $m=10000$, $\ell=2500$, $n=5000$, and $s=20$. The relative residual $\|\mbf b-\mbf A\mbf B\mbf x^{(k)}\|_2/\|\mbf b\|_2$, the relative error $\|\mbf x^{(k)}-\mbf x_\star\|_2/\|\mbf x_\star\|_2$,  and approximated solutions (last iterates of RK-RK and RK-RSK) are averaged over 50 independent runs.  We have the following observations: (i) the RK-RK algorithm converges to a solution, but not a sparse one; (ii) the RK-RSK algorithm converges to a sparse solution and indeed recovers $\mbf x_\star$; (iii) the RK-RSK algorithm has a better convergence rate compared with the RK-RK algorithm. 
 
In Figure \ref{fig2}, we plot the results for the case $\mbf b\notin\ran(\mbf A\mbf B)$ with $m=10000$, $\ell=2500$, $n=5000$, and $s=20$. The relative residual $\|\mbf B^\top\mbf A^\top(\mbf b-\mbf A\mbf B\mbf x^{(k)})\|_2/\|\mbf B^\top\mbf A^\top\mbf b\|_2$, the relative error $\|\mbf x^{(k)}-\mbf x_\star\|_2/\|\mbf x_\star\|_2$, and approximated solutions (last iterates of RGS-RK and RGS-RSK) are averaged over 50 independent runs. We have the following observations: (i) the RGS-RK algorithm converges to a least squares solution, but not a sparse one; (ii) the RGS-RSK algorithm converges to a sparse least squares solution and indeed recovers $\mbf x_\star$; (iii) the RGS-RSK algorithm has a better convergence rate compared with the RGS-RK algorithm.

\subsection{Example 2} We consider the wine quality data set obtained from the UCI Machine Learning Repository \cite{dua2017uci}. Let $\mbf X\in\mbbr^{m\times n}$ denote the data matrix (a sample of $m=1599$ red wines with $n=11$ physio-chemical properties of each wine). The matrices $\mbf A\in\mbbr^{m\times 5}$ and $\mbf B\in\mbbr^{5\times n}$ are obtained by using the MATLAB's nonnegative matrix factorization function {\tt nnmf} as follows: {\tt [A,B]=nnmf(X,5)}. We compute $\mbf C=\mbf A\mbf B$ in MATLAB. Let $\mbf x_\star\in\mbbr^{11}$ be a 3-sparse vector with support $\{1,6,11\}$. The three nonzero entries of $\mbf x_\star$ are set to be 1. The vector $\mbf b$ is constructed by using the same approach as that of Example 1. We compare the RK-RSK algorithm and the RGS-RSK algorithm for the factorized linear system $\mbf A\mbf B\mbf x=\mbf b$ with the RSK algorithm and the GERK-(a,d) algorithm for the full linear system $\mbf C\mbf x=\mbf b$. For the proposed algorithms, we use $\lambda=1$, $\mbf y^{(0)}=\mbf 0$, $\mbf z^{(0)} = \mbf 0$, and the maximum number of iterations {\tt maxit=10$m$}.
 
\begin{figure}[htb]
\centerline{\epsfig{figure=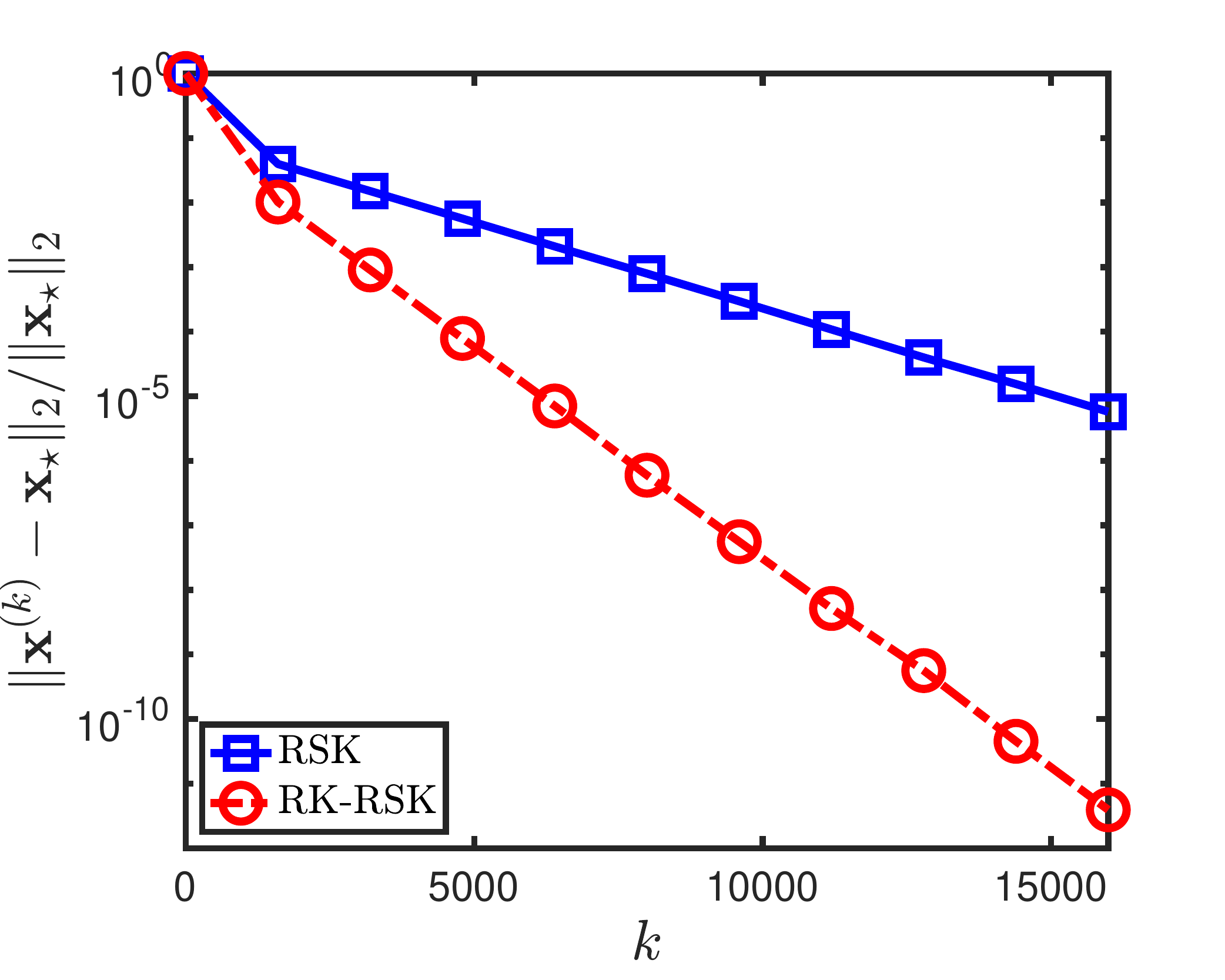,height=1.65in}\epsfig{figure=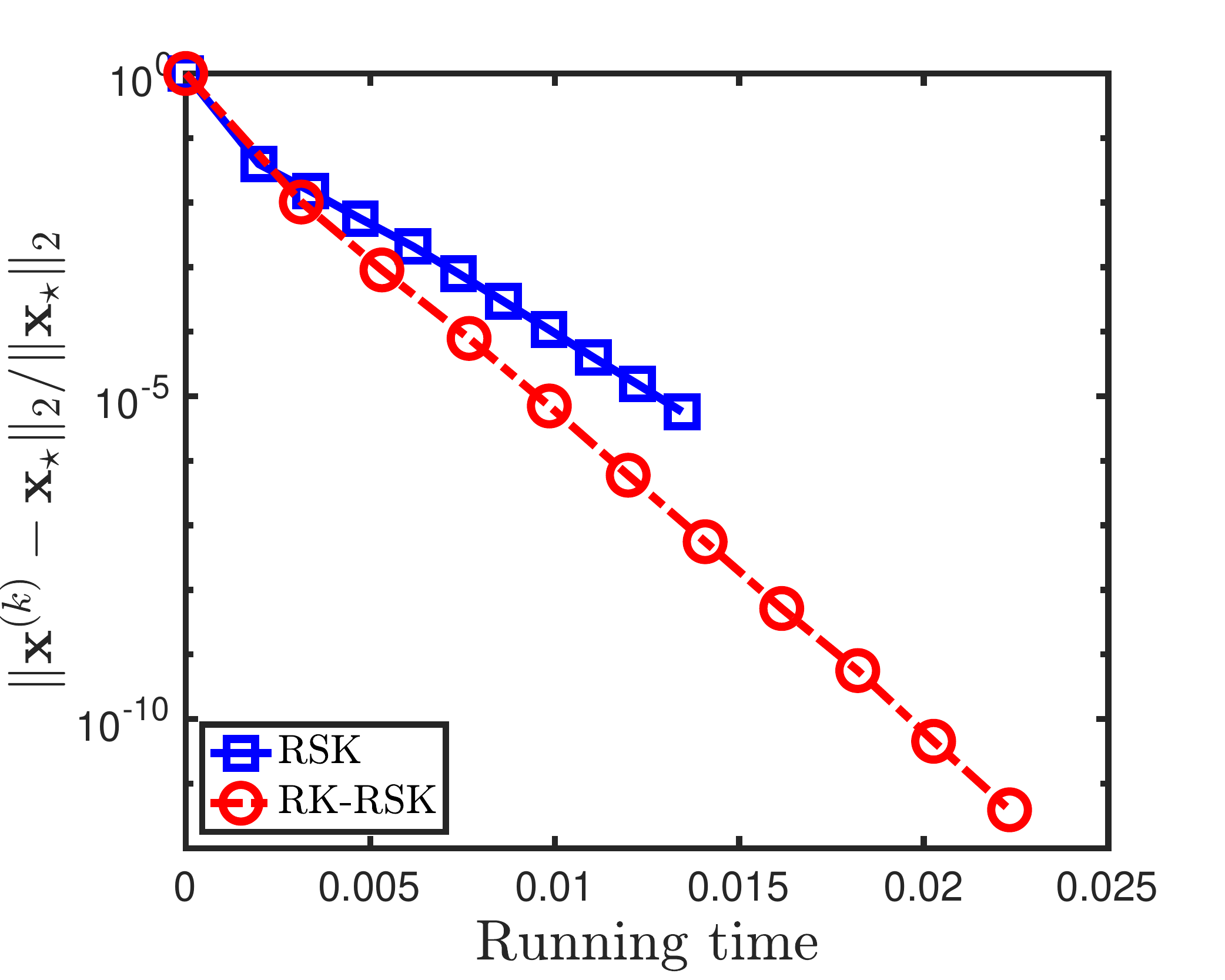,height=1.65in}\epsfig{figure=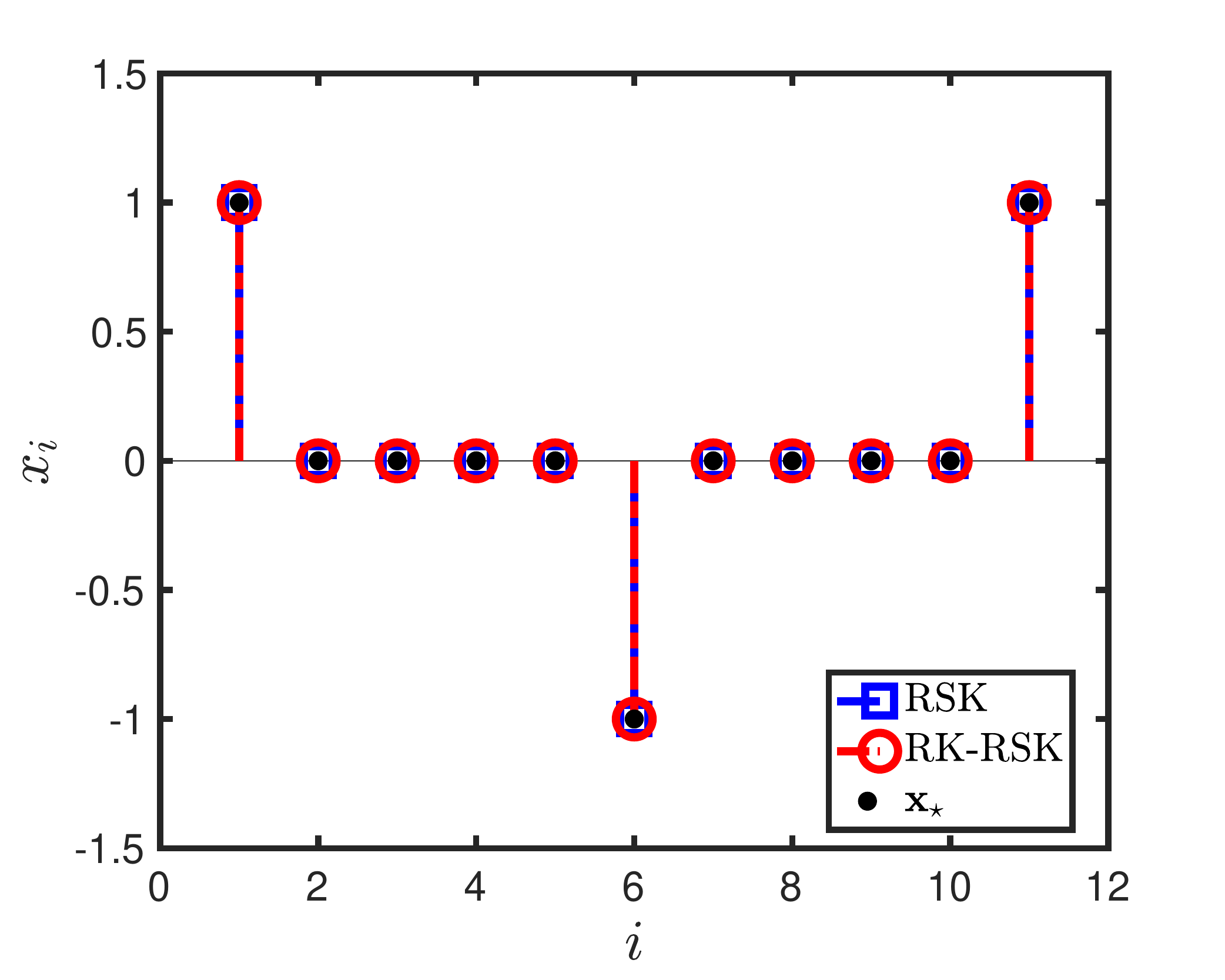,height=1.65in}}
\caption{Comparison of RSK and RK-RSK for the wine quality data set. Left: iteration vs relative error $\|\mbf x^{(k)}-\mbf x_\star\|_2/\|\mbf x_\star\|_2$.  Middle: running time (in seconds) vs relative error $\|\mbf x^{(k)}-\mbf x_\star\|_2/\|\mbf x_\star\|_2$. Right: approximated solutions (last iterates of RSK and RK-RSK) and the sparse vector $\mbf x_\star$.}
\label{fig3}
\end{figure}  

\begin{figure}[htb]
\centerline{\epsfig{figure=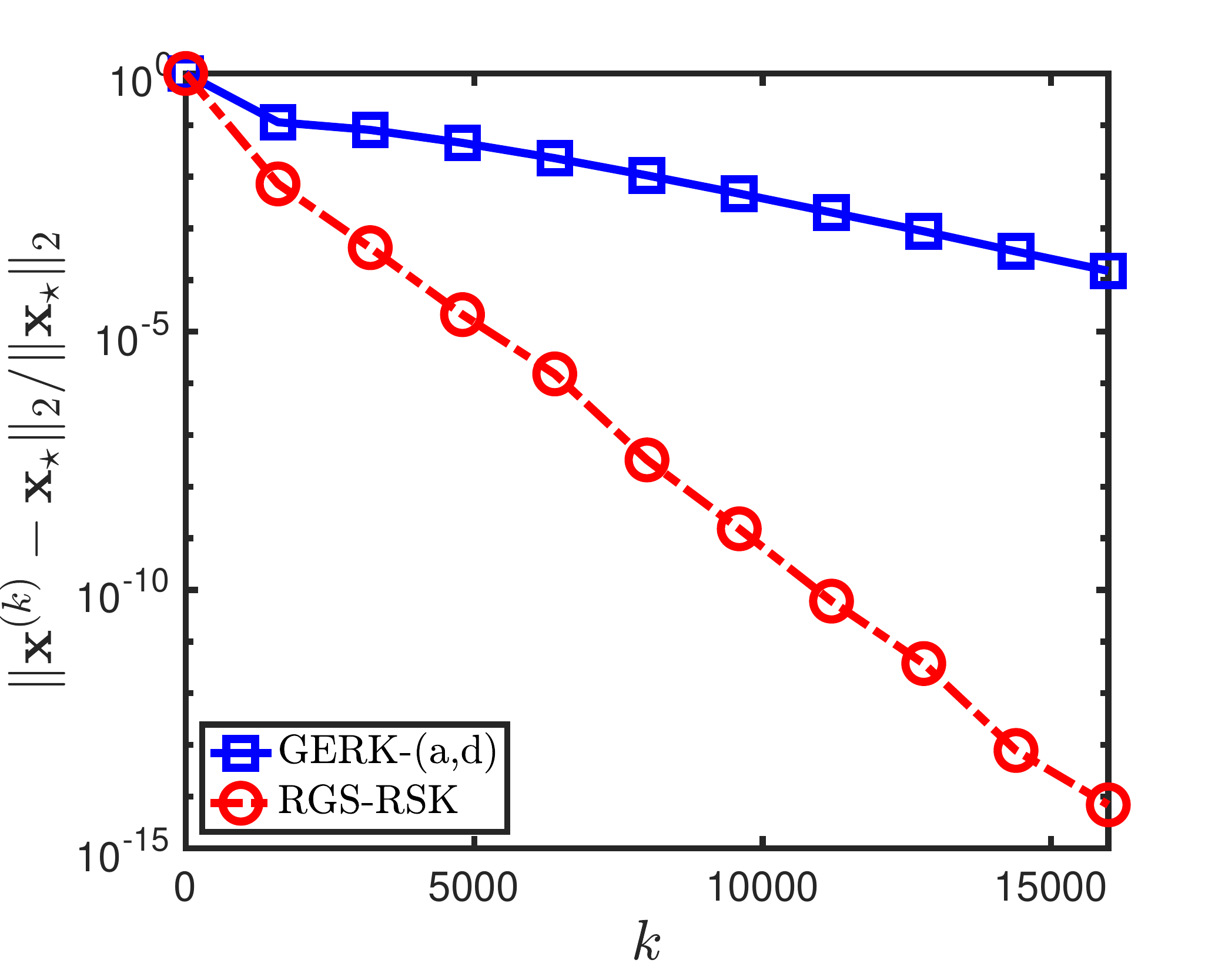,height=1.65in}\epsfig{figure=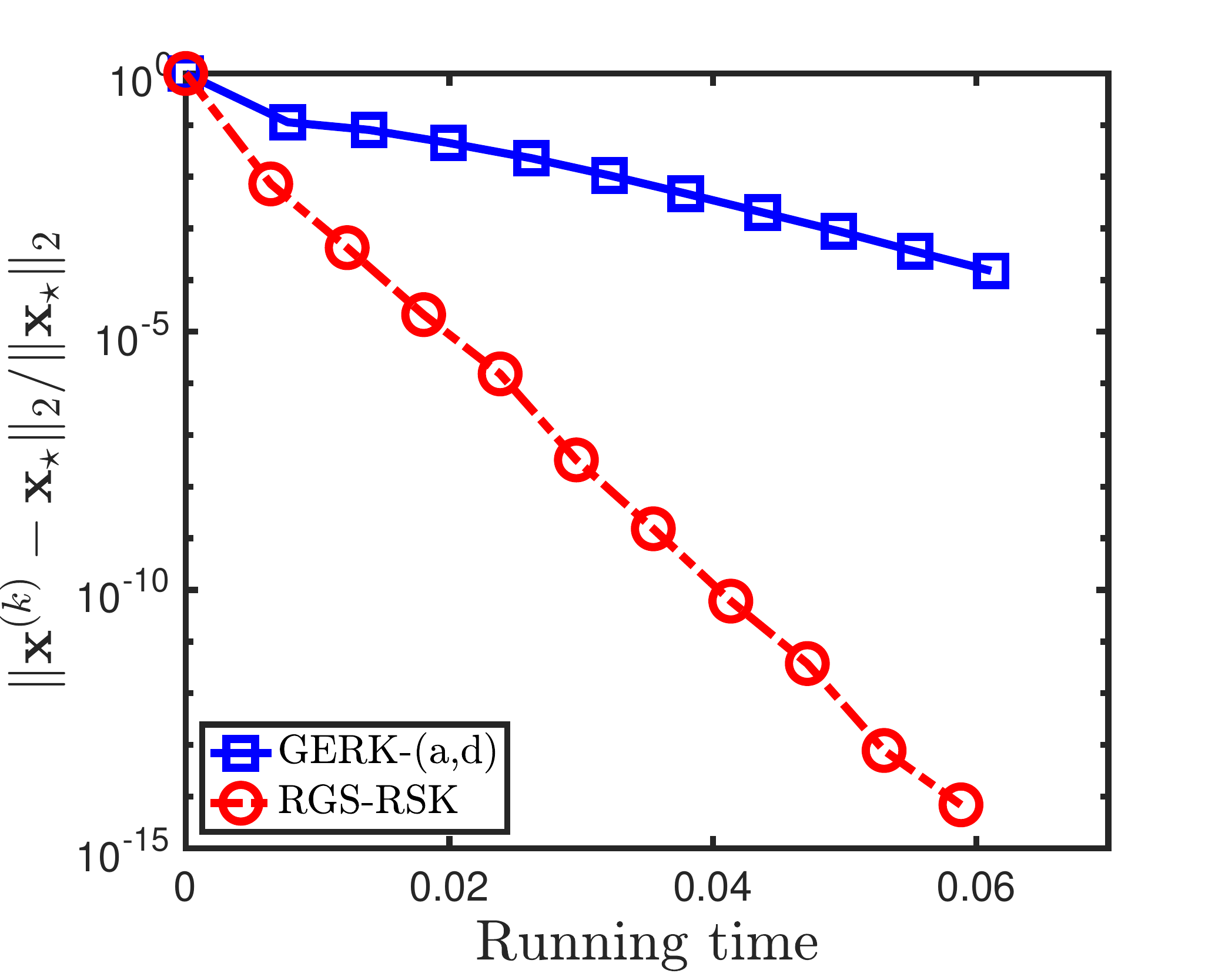,height=1.65in}\epsfig{figure=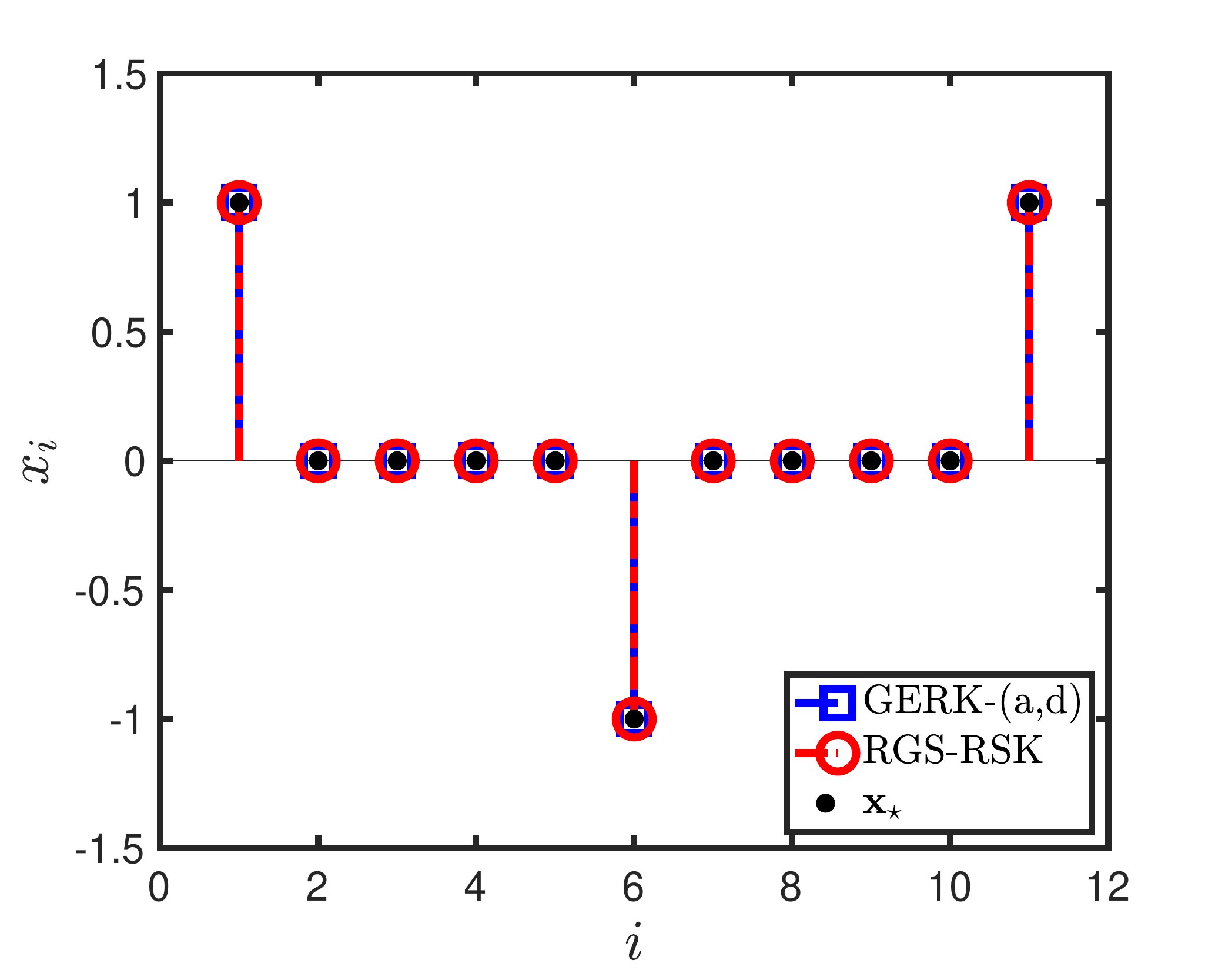,height=1.65in}}
\caption{Comparison of GERK-(a,d) and RGS-RSK for the wine quality data set. Left: iteration vs relative error $\|\mbf x^{(k)}-\mbf x_\star\|_2/\|\mbf x_\star\|_2$.  Middle: running time (in seconds) vs relative error $\|\mbf x^{(k)}-\mbf x_\star\|_2/\|\mbf x_\star\|_2$. Right: approximated solutions (last iterates of GERK-(a,d) and RGS-RSK) and the sparse vector $\mbf x_\star$.}
\label{fig4}
\end{figure} 

In Figures \ref{fig3}, we plot iteration vs relative error $\|\mbf x^{(k)}-\mbf x_\star\|_2/\|\mbf x_\star\|$, running time (in seconds) vs relative error $\|\mbf x^{(k)}-\mbf x_\star\|_2/\|\mbf x_\star\|$, and  approximated solutions (last iterates of RSK and RK-RSK). The results are averaged over 50 independent runs. We have the following observations: (i) both the RSK algorithm and the RK-RSK algorithm recover the sparse solution $\mbf x_\star$; (ii) the RK-RSK algorithm is faster than the RSK algorithm. Figure \ref{fig4} reports the results for the GERK-(a,d) algorithm and the RGS-RSK algorithm. And we have the following observations: (i) both the GERK-(a,d) algorithm and the RGS-RSK algorithm recover the sparse least squares solution $\mbf x_\star$; (ii) the RGS-RSK algorithm is faster than the GERK-(a,d) algorithm. 
 
\section{Concluding remarks} We have proposed two algorithms to find solutions with certain structures of a factorized linear system. We have proved their linear convergence under some assumptions. Our numerical examples indicate that the proposed algorithms can find sparse (least squares) solutions and can be faster than the RSK and GERK-(a,d) algorithms for the corresponding full linear system. Existing acceleration strategies for the RK, RGS, and RRK algorithms such as those used in \cite{needell2014paved,bai2018greed,bai2019greed,necoara2019faste,zhang2022greed,jiang2022rando,tondji2022faste,yuan2022spars,yuan2022adapt,wu2022two} can be integrated into our algorithms easily and the corresponding convergence analysis is straightforward. The extension to a factorized linear system  with rank-deficient $\mbf A$ and $\mbf B$ will be the future work. 
\
\section*{Acknowledgments} The author thanks the referees for detailed comments and suggestions that have led to significant improvements. The author also thanks Xuemei Chen and Jing Qin for their code. This work was supported by the National Natural Science Foundation of China (No.12171403 and No.11771364), the Natural Science Foundation of Fujian Province of China (No.2020J01030), and the Fundamental Research Funds for the Central Universities (No.20720210032).
 
 \section*{Data availability}
The data that support the findings of this study are available from the UCI Machine Learning Repository, http://archive.ics.uci.edu/ml.

\section*{Conflict of Interest}
The author declares that he has no conflict of interest.

{\small 
%\bibliographystyle{plain}
%\bibliography{RRIM.bib}

}

\end{document}